\documentclass[11pt,a4paper]{article}

\usepackage[a4paper,margin=2.5cm]{geometry}
\usepackage[T1]{fontenc}
\usepackage[utf8]{inputenc}
\usepackage{lmodern}
\usepackage{microtype}
\usepackage{longtable}

\usepackage{amssymb,amsmath,amsfonts,latexsym}
\usepackage{amsthm}
\usepackage{mathtools}
\usepackage{bm}

\usepackage{graphicx}
\usepackage{subcaption}
\usepackage{float}
\usepackage{multirow,array}
\usepackage{booktabs}
\usepackage{xcolor}
\usepackage{tikz}
\usetikzlibrary{arrows.meta,positioning,calc}

\usepackage{enumitem}
\allowdisplaybreaks

\usepackage[colorlinks=true,
linkcolor=blue,
citecolor=blue,
urlcolor=blue,
bookmarksnumbered=true]{hyperref}

\numberwithin{equation}{section}

\newtheorem{theorem}{Theorem}[section]
\newtheorem{lemma}[theorem]{Lemma}
\newtheorem{proposition}[theorem]{Proposition}
\newtheorem{corollary}[theorem]{Corollary}
\newtheorem{definition}[theorem]{Definition}
\newtheorem{example}[theorem]{Example}
\newtheorem{remark}[theorem]{Remark}

\title{Characterizations and Algebraic, Perturbation, and Topological Properties of $\theta(H)$}

\author{Zinan Su$^{*}$ \and Yuanhang Zhang$^{\dagger}$ \\[1ex]
	\small School of Mathematics, Jilin University, Qianjin Street 2699,
	Changchun 130012, Jilin, China \\
	\small $^*$Corresponding author: \texttt{suzinan20050803@aliyun.com}\\
\small $\dagger$: \texttt{zhangyuanhang@jlu.edu.cn}}

\date{}

\begin{document}
	
	\maketitle
	
	\begin{abstract}
		Let $H$ be a separable complex Hilbert space, and let $\theta(H)$ be the
		class of bounded operators whose restrictions to every nonzero closed
		reducing subspace attain their numerical radii. We characterize this
		property by eigenvalues of suitable self-adjoint parts. For normal
		operators, membership in $\theta(H)$ is equivalent to pure point spectrum
		and reverse well-ordering of the moduli of eigenvalues. Without normality
		assumptions, we prove stability under finite orthogonal direct sums using
		graph subspaces of closed, possibly unbounded, intertwining operators.
		We give a compact-perturbation criterion in terms of the essential
		numerical radius and show that an operator remains in $\theta(H)$ under every compact perturbation if and only if it is compact. Finally, $\theta(H)$ is norm-dense in $\mathcal B(H)$ and is norm closed if and only
		if $H$ is finite-dimensional.
	\end{abstract}
	
	
	
	\setcounter{tocdepth}{1}
	\tableofcontents
	\clearpage

	\section{Introduction}\label{sec:intro}
	
	Attainment of extremal values is a recurring theme in functional analysis.
	Even when an operator-theoretic quantity is defined by a supremum over
	unit vectors, there need not be a vector at which that supremum is achieved.
	For the operator norm, this question belongs to the classical theory of
	norm-attaining functionals and operators, initiated by the Bishop--Phelps
	theorem and subsequently developed in the operator setting
	\cite{bishop1961,Lindenstrauss1963,Carvajal2012}. A related, but genuinely different,
	problem arises for the numerical radius. The distinction is particularly
	important when one asks whether attainment persists after passing from an
	operator to the pieces singled out by its invariant structure.
	
	Let $H$ be a nonzero separable complex Hilbert space, with inner products
	linear in the first variable, and write $\mathcal B(H)$ for the
	bounded linear operators on $H$. For $T\in\mathcal B(H)$, the numerical
	range and numerical radius are
	\[
	W(T)=\{\langle Tx,x\rangle:\|x\|=1\},
	\qquad
	\omega(T)=\sup_{\|x\|=1}|\langle Tx,x\rangle|.
	\]
	The numerical range is a classical geometric invariant of an operator:
	its convexity is the Toeplitz--Hausdorff theorem
	\cite{toeplitz1918,hausdorff1919}, and its relationship with spectral
	and norm properties is developed in
	\cite{bonsall1971,bonsall1973,gustafson1997}. In particular,
	$\omega$ is a norm equivalent to the operator norm, since
	\[
	\frac12\|T\|\leqslant\omega(T)\leqslant\|T\|.
	\]
	These estimates control the \emph{size} of the numerical radius,
	but do not decide whether its defining supremum is attained.
	Indeed, for $\zeta$ on the unit circle $\mathbb{T}$, consider the
	self-adjoint part
	\[
	H_\zeta(T)=\frac{\zeta T+\overline{\zeta}T^*}{2}.
	\]
	The identity
	\[
	\omega(T)=\max_{\zeta\in\mathbb{T}}\sup\sigma(H_\zeta(T))
	\]
	shows that numerical-radius attainment is a boundary spectral
	question: it is equivalent to the existence of a suitable
	$H_\zeta(T)$ having $\omega(T)$ as an eigenvalue
	(Lemma~\ref{lem:numerical}). Thus the issue is not merely the value
	of an extremum, but whether its supporting spectral boundary
	contains an eigenvector.
	
	This distinction becomes visible in infinite dimensions.
	For example, the positive diagonal operator
	$De_n=\left(1-\frac1n\right)e_n$ on $\ell^2(\mathbb N)$ satisfies
	$\omega(D)=1$, although no unit vector attains this value.
	More importantly, even \emph{global} numerical-radius attainment
	does not ensure attainment on a part of the operator:
	$1\oplus D$ attains its numerical radius on
	$\mathbb C\oplus\ell^2(\mathbb N)$, whereas its restriction to
	the reducing summand $\{0\}\oplus\ell^2(\mathbb N)$ does not.
	This example leads to a stronger question:
	\emph{which operators attain their numerical radii not only globally,
		but on every nonzero closed reducing subspace?}
	
	For $T\in\mathcal B(H)$, let $\operatorname{Red}(T)$ denote the
	family of its closed reducing subspaces. We study the class
	\[
	\theta(H)=\left\{T\in\mathcal B(H):
	\begin{array}{l}
		\text{for every }M\in\operatorname{Red}(T)\setminus\{0\},
		\text{ there exists }x\in M,\\
		\|x\|=1,\quad
		|\langle Tx,x\rangle|=\omega(T|_M)
	\end{array}\right\}.
	\]
	The requirement is naturally tied to operator structure.
	A closed subspace $M$ reduces $T$ precisely when its orthogonal
	projection $P_M$ commutes with both $T$ and $T^*$.
	Consequently, the reducing-subspace lattice is the projection
	lattice of the von Neumann algebra $\{T,T^*\}'$.
	Membership in $\theta(H)$ can therefore be viewed as a
	\emph{hereditary attainment condition} indexed by the projections
	in this commutant. The use of reducing subspaces is essential:
	unlike an arbitrary compression, each restriction $T|_M$
	is an operator on a genuine orthogonal summand, with
	$(T|_M)^*=T^*|_M$. The problem consequently connects
	the geometry of the numerical range with the decomposition
	theory of operators.
	
	The first issue is spectral. For a normal operator $N$,
	$\omega(N)=\|N\|$, but numerical-radius attainment is still
	more restrictive than norm attainment. For instance, the
	bilateral shift is unitary and hence attains its operator norm,
	yet it has no eigenvalues and does not attain its numerical
	radius. When attainment is required on \emph{all} reducing
	subspaces, the spectral restrictions are stronger still.
	Our normal-operator characterization
	(Theorem~\ref{thm:normal}) identifies two precise conditions:
	the eigenspaces must have closed linear span $H$, and the set
	\[
	\{|\lambda|:\lambda\in\sigma_p(N)\}
	\]
	must be reverse well-ordered, meaning that every nonempty
	subset has a largest element. The first condition excludes
	a residual continuous spectral part; the second prevents a
	reducing spectral summand from having an unattained supremum
	of eigenvalue moduli. In this way, a geometric extremal
	property becomes a condition on both spectral type and the
	order structure of the point spectrum.
	
	A second issue concerns algebraic stability. If an attainment
	property is imposed on every reducing subspace, it is not
	obvious that it should survive the formation of direct sums.
	Even for $T_1\oplus T_2$, a reducing subspace need not be a
	coordinatewise sum of reducing subspaces of $T_1$ and $T_2$:
	the two coordinates can be coupled by a graph subspace.
	Such graphs may be associated with closed \emph{unbounded}
	intertwining operators, so an argument confined to bounded
	graph maps would miss legitimate reducing subspaces.
	The structure of reducing subspaces of direct sums is itself
	a significant topic \cite{FangWuJiang2003}. Nevertheless,
	we prove that $\theta$ is preserved by finite orthogonal
	direct sums without assuming normality
	(Theorem~\ref{thm:finite-direct-sum}). This contrast reflects the fact that, for infinitely many summands, the supremum of the component numerical radii need not be attained.
	
	A third issue is perturbation theory. The essential numerical
	range, introduced by Fillmore, Stampfli and Williams
	\cite{FillmoreStampfliWilliams1972}, records numerical-range
	behavior that persists under compact perturbations. Its
	associated essential numerical radius $\omega_e(T)$ provides
	a natural threshold: when $\omega(T)>\omega_e(T)$, the
	numerical radius is attained, because the relevant
	self-adjoint part has an eigenvalue above its essential
	spectral boundary. For the class $\theta(H)$, however, one
	must examine this comparison on every reducing subspace
	of the \emph{perturbed} operator, and the reducing subspaces
	themselves may change under perturbation. This makes it
	meaningful to ask not just whether individual compact
	perturbations preserve membership, but which operators
	remain in $\theta(H)$ under \emph{every} compact perturbation.
	We show that this compact-perturbation-invariant core is
	exactly $\mathcal K(H)$ (Theorem~\ref{thm:invariant}),
	and we exhibit a non-normal irreducible member of
	$\theta(H)$ whose membership is lost after a rank-one
	perturbation (Proposition~\ref{prop:rankone-instability}).
	These results distinguish an attainment property from
	robustness of that property.
	
	Finally, there is an approximation question. The
	Bishop--Phelps tradition asks how large an attainment class
	is in the ambient norm topology; for numerical-radius
	attainment, Berg and Sims established an important density
	result \cite{BergSims1984}. Requiring attainment on every
	reducing subspace appears to impose a much stronger
	condition. Yet irreducibility changes the picture:
	for an irreducible operator, there is no proper nonzero
	reducing subspace, so the hereditary requirement reduces
	to ordinary numerical-radius attainment. Combining this
	observation with Halmos's density theorem for irreducible
	operators \cite{Halmos1968} and an attainment argument,
	we prove that $\theta(H)$ is norm-dense in $\mathcal B(H)$
	(Theorem~\ref{thm:norm-density}). On an infinite-dimensional
	space it is nevertheless not norm closed, and its
	compact-perturbation-invariant core is much smaller.
	Thus \emph{density}, \emph{hereditary attainment}, and
	\emph{perturbation stability} are three distinct aspects
	of the same problem.
	
	Section~\ref{sec:general} gives the general attainment criterion.
	Sections~\ref{sec:normal}--\ref{sec:finite-direct-sum} develop the normal
	characterization and algebraic properties, including finite direct sums.
	Section~\ref{sec:compact} treats compact perturbations and their invariant
	core, and Section~\ref{sec:topology} establishes the topological results.
	Appendix~\ref{app:unbounded} collects the technical closed-operator
	lemmas used in the finite direct-sum argument.
	
	\section{Preliminaries and a characterization of general operators}\label{sec:general}

	Throughout the paper, $H$ is a nonzero separable complex Hilbert space, not necessarily infinite-dimensional. We use $\mathbb N=\{1,2,\ldots\}$. We write $\mathcal{B}(H)$ for the algebra of all bounded linear operators on $H$ and $\mathcal{K}(H)$ for the two-sided ideal of compact operators in $\mathcal{B}(H)$. When $H$ is finite-dimensional, $\mathcal{K}(H)=\mathcal{B}(H)$ and every operator belongs to $\theta(H)$. For a closed subspace $M$ of $H$, let $P_M$ be the orthogonal projection onto $M$. A closed subspace $M$ is said to reduce $T\in\mathcal{B}(H)$ if
	\[
	T(M)\subset M
	\qquad\text{and}\qquad
	T^*(M)\subset M.
	\]
	We call $T$ \emph{irreducible} if its only closed reducing subspaces are
	$\{0\}$ and $H$. Equivalently, the commutant of $T$ and $T^*$ satisfies
	$\{T,T^*\}'=\mathbb{C}I$. For an operator $T\in\mathcal{B}(H)$, we denote by $\sigma(T)$, $\sigma_p(T)$ and $\sigma_e(T)$ its spectrum, its point spectrum and its Fredholm essential spectrum, respectively. We shall refer to $\sigma_e(T)$ as the essential spectrum.

	We shall use the classical numerical-radius bounds
	\begin{equation}\label{eq:numerical-radius-bounds}
		\frac12\|T\|\leqslant\omega(T)\leqslant\|T\|,
		\qquad T\in\mathcal B(H);
	\end{equation}
	see \cite{bonsall1971,dragomir2013}. Both constants are sharp: for
	$J=\begin{pmatrix}0&1\\0&0\end{pmatrix}$ on $\mathbb C^2$,
	$\|J\|=1$ and $\omega(J)=\frac12$, whereas $\omega(I)=\|I\|=1$.
	For refinements of numerical-radius estimates, see
	\cite{kittaneh2003,yamazaki2007}.
	
	We shall repeatedly use the following standard facts:
	\begin{itemize}
		\item[(F1)] If $A$ is self-adjoint, then
		\[
		\sup\sigma(A)=\sup_{\|x\|=1}\langle Ax,x\rangle.
		\]
		\item[(F2)] If $N$ is normal, then
		\[
		\omega(N)=\|N\|.
		\]
	\end{itemize}
	
	\begin{definition}\label{def:rwo}
		A totally ordered set $(X,<)$ is called \emph{reverse well-ordered} if every nonempty subset of $X$ has a maximum element. Equivalently, the reverse order is a well-order. For $X\subset\mathbb{R}$, this is equivalent to the nonexistence of a strictly increasing infinite sequence in $X$ (the equivalence for subsets of $\mathbb{R}$ uses a routine countable choice argument).
	\end{definition}
	
	\begin{definition}\label{def:purepoint}
		An operator $T\in\mathcal{B}(H)$ is said to have \emph{pure point spectrum} if
		\[
		H=
		\overline{\operatorname{span}}
		\left(
		\bigcup_{\lambda\in\sigma_p(T)}
		\ker(T-\lambda I)
		\right).
		\]
		When $T$ is normal, this is equivalent to the existence of an orthonormal basis of eigenvectors of $T$.
	\end{definition}
	
	Before introducing the self-adjoint operators $H_\zeta(S)$ and proving the
	general characterization, we record a basic fact that will be used repeatedly:
	every compact operator belongs to $\theta(H)$. Its proof is a direct compactness
	argument and does not rely on Theorem~\ref{thm:general}.
	
	\begin{proposition}\label{prop:compact}
		Every compact operator belongs to $\theta(H)$.
	\end{proposition}
	
	\begin{proof}
		Let $T\in\mathcal K(H)$, take $0\neq M\in\operatorname{Red}(T)$, and put $S=T|_M$. Since $T$ is compact and $M$ is closed and reducing, $S$ is compact. Let $r=\omega(S)$. If $r=0$, then $S=0$ because the numerical radius is a norm, and the attainment on $M$ is immediate.
		
		Assume $r>0$. Choose a sequence of unit vectors $\{x_n\}\subset M$ such that
		\[
		|\langle Sx_n,x_n\rangle|\to r.
		\]
		By weak compactness of the unit ball, after passing to a subsequence we may assume that $x_n\rightharpoonup x$, where $x_n\rightharpoonup x$ means that the sequence $\{x_n\}$ converges weakly to $x$. Since $M$ is weakly closed, $x\in M$, and compactness of $S$ gives $Sx_n\to Sx$ in norm. Therefore
		\[
		\langle Sx_n,x_n\rangle-\langle Sx,x\rangle
		=
		\langle Sx_n-Sx,x_n\rangle
		+
		\langle Sx,x_n-x\rangle
		\to0.
		\]
		Hence $|\langle Sx,x\rangle|=r$. Since
		\[
		|\langle Sy,y\rangle|\leqslant r\|y\|^2
		\qquad (y\in M),
		\]
		we obtain
		\[
		r=|\langle Sx,x\rangle|\leqslant r\|x\|^2.
		\]
		Thus $\|x\|\ge1$, while weak lower semicontinuity gives
		\[
		\|x\|\leqslant\varliminf_{n\to\infty}\|x_n\|=1.
		\]
		Consequently $\|x\|=1$, so $\omega(S)$ is attained on $M$. Since $M$ was arbitrary, $T\in\theta(H)$.
	\end{proof}
	
	For $S\in\mathcal{B}(H)$ and $\zeta\in\mathbb{T}$, define
	\[
	H_\zeta(S)=\frac{\zeta S+\overline{\zeta}S^*}{2},
	\]
	where $\mathbb{T}$ denotes the unit circle in $\mathbb{C}$. Clearly, $H_\zeta(S)$ is self-adjoint. The following lemma, which is the key tool of this section, expresses the numerical radius and the numerical-radius attainment of an arbitrary bounded operator through the top of the spectrum of these self-adjoint operators.
	
	\begin{lemma}\label{lem:numerical}
		For every $S\in\mathcal{B}(H)$,
		\[
		\omega(S)=\max_{\zeta\in\mathbb{T}}\sup\sigma(H_\zeta(S)).
		\]
		Moreover, $\omega(S)$ is attained on the unit sphere if and only if there exist $\zeta\in\mathbb{T}$ and a unit vector $x\in H$ such that
		\[
		H_\zeta(S)x=\omega(S)x.
		\]
	\end{lemma}
	
	\begin{proof}
		If $\omega(S)=0$, then $S=0$ because the numerical radius is a norm.
		Both assertions are immediate, since $H_\zeta(S)=0$ for every
		$\zeta\in\mathbb{T}$. Assume henceforth that $\omega(S)>0$.
		
		For every unit vector $x$,
		\[
		\sup_{\zeta\in\mathbb{T}}
		\langle H_\zeta(S)x,x\rangle
		=
		\sup_{\zeta\in\mathbb{T}}
		\operatorname{Re}\left(\zeta\langle Sx,x\rangle\right)
		=
		|\langle Sx,x\rangle|.
		\]
		Hence
		\[
		\begin{aligned}
			\omega(S)=\sup_{\|x\|=1}|\langle Sx,x\rangle|=\sup_{\|x\|=1}\sup_{\zeta\in\mathbb{T}}
			\langle H_\zeta(S)x,x\rangle=\sup_{\zeta\in\mathbb{T}}\sup_{\|x\|=1}
			\langle H_\zeta(S)x,x\rangle.
		\end{aligned}
		\]
		By (F1),
		\[
		\omega(S)=
		\sup_{\zeta\in\mathbb{T}}\sup\sigma(H_\zeta(S)).
		\]
		For each $\zeta$,
		\[
		\sup\sigma(H_\zeta(S))
		=
		\sup_{z\in\overline{W(S)}}\operatorname{Re}(\zeta z).
		\]
		The right-hand side is the support function of the compact convex set $\overline{W(S)}$, hence it is continuous in $\zeta$. Since $\mathbb{T}$ is compact, the supremum is attained.
		
		Suppose first that $\omega(S)$ is attained at a unit vector $x$. Choose $\zeta\in\mathbb{T}$ such that
		\[
		\zeta\langle Sx,x\rangle
		=
		|\langle Sx,x\rangle|
		=
		\omega(S).
		\]
		Then
		\[
		\langle H_\zeta(S)x,x\rangle=\omega(S).
		\]
		For every unit vector $y$,
		\[
		\langle H_\zeta(S)y,y\rangle
		\leqslant|\langle Sy,y\rangle|
		\leqslant\omega(S).
		\]
		Thus $H_\zeta(S)\leqslant\omega(S)I$, and equality at $x$ implies
		\[
		H_\zeta(S)x=\omega(S)x.
		\]
		Conversely, if $H_\zeta(S)x=\omega(S)x$ for some unit vector $x$, then
		\[
		\omega(S)
		=
		\operatorname{Re}\left(\zeta\langle Sx,x\rangle\right)
		\leqslant
		|\langle Sx,x\rangle|
		\leqslant\omega(S),
		\]
		so $|\langle Sx,x\rangle|=\omega(S)$.
	\end{proof}
	
	By applying Lemma~\ref{lem:numerical} to the restriction of $T$ to each reducing subspace, we obtain the following characterization of the class $\theta(H)$ for arbitrary bounded operators.
	
	\begin{theorem}\label{thm:general}
		Let $T\in\mathcal{B}(H)$. Then $T\in\theta(H)$ if and only if, for every nonzero $M\in\operatorname{Red}(T)$, there exists $\zeta\in\mathbb{T}$ such that
		\[
		\omega(T|_M)\in
		\sigma_p\left(H_\zeta(T|_M)\right).
		\]
		Equivalently, for every nonzero $M\in\operatorname{Red}(T)$, there exist
		$\zeta\in\mathbb{T}$ and a unit vector $x\in M$ such that
		\[
		H_\zeta(T|_M)x=\omega(T|_M)x.
		\]
	\end{theorem}
	
	\begin{proof}
		By the definition of $\theta(H)$, the assertion $T\in\theta(H)$ is equivalent to the attainment of $\omega(T|_M)$ on the unit sphere of every nonzero reducing subspace $M$. Applying Lemma~\ref{lem:numerical} to $S=T|_M$ gives the result.
	\end{proof}

	\section{A spectral characterization for normal operators}\label{sec:normal}

	In this section we restrict attention to normal operators and obtain a complete spectral description of membership in $\theta(H)$. The main result, Theorem~\ref{thm:normal}, states that a normal operator belongs to $\theta(H)$ if and only if it has pure point spectrum and the set of moduli of its eigenvalues is reverse well-ordered. The proof proceeds in two steps: first we show that membership in $\theta(H)$ forces pure point spectrum, and then we show that, under pure point spectrum, the reverse well-order condition on the eigenvalue moduli is both necessary and sufficient.
	
	\begin{lemma}\label{lem:purepoint}
		Let $T$ be normal. If $T\in\theta(H)$, then $T$ has pure point spectrum.
	\end{lemma}
	
	\begin{proof}
		Let
		\[
		K=
		\overline{\operatorname{span}}
		\{\ker(T-\lambda I):\lambda\in\sigma_p(T)\}.
		\]
		For a normal operator, every eigenspace $\ker(T-\lambda I)$ reduces $T$, and hence $K$ also reduces $T$. Suppose that $T$ does not have pure point spectrum. Then $K\neq H$. Put
		\[
		M=K^\perp.
		\]
		Then $M\neq\{0\}$ and $M\in\operatorname{Red}(T)$. Moreover, $T|_M$ has no eigenvalues. Indeed, if $T|_M x=\lambda x$ for some nonzero $x\in M$, then $x\in\ker(T-\lambda I)\subset K$, contradicting $x\in K^\perp$.
		
		Since $T|_M$ is normal, (F2) gives
		\[
		\omega(T|_M)=\|T|_M\|.
		\]
		As $T\in\theta(H)$, there is a unit vector $x\in M$ such that
		\[
		|\langle Tx,x\rangle|=\|T|_M\|.
		\]
		By Cauchy--Schwarz,
		\[
		|\langle Tx,x\rangle|
		\leqslant \|Tx\|\|x\|
		\leqslant \|T|_M\|.
		\]
		Equality throughout implies
		\[
		|\langle Tx,x\rangle|=\|Tx\|\|x\|.
		\]
		The equality condition in Cauchy--Schwarz yields $Tx=\lambda x$ for some $\lambda\in\mathbb{C}$, contradicting the fact that $T|_M$ has no eigenvalues. Therefore $T$ has pure point spectrum.
	\end{proof}
	
	The next lemma shows that, for a normal operator with pure point spectrum, every reducing subspace decomposes along the eigenspaces of the operator. This is the second structural ingredient in the proof of Theorem~\ref{thm:normal}.
	
	\begin{lemma}\label{lem:reducing}
		Let $T$ be normal with pure point spectrum. For every nonzero $M\in\operatorname{Red}(T)$, put
		\[
		M_\lambda=M\cap\ker(T-\lambda I),
		\qquad
		\lambda\in\sigma_p(T).
		\]
		Then
		\[
		M=
		\bigoplus_{\lambda\in\sigma_p(T)}M_\lambda.
		\]
	\end{lemma}
	
	\begin{proof}
		For each $\lambda\in\sigma_p(T)$, let $P_\lambda$ be the orthogonal projection onto
		\[
		E_\lambda=\ker(T-\lambda I).
		\]
		Since $T$ is normal, $P_\lambda$ is a spectral projection. If $P_M$
		denotes the projection onto $M$, then $M\in\operatorname{Red}(T)$ implies
		that $P_M$ commutes with $T$ and $T^*$, and hence with every spectral
		projection of $T$. In particular, $P_M$ commutes with $P_\lambda$, so
		\[
		P_\lambda(M)=P_\lambda P_M(H)
		=P_M P_\lambda(H)
		\subset M.
		\]
		Since $P_\lambda(H)=E_\lambda$, we have
		\[
		P_\lambda(M)\subset M\cap E_\lambda=M_\lambda.
		\]
		Conversely, if $x\in M_\lambda$, then $P_\lambda x=x$, so $x\in P_\lambda(M)$. Thus
		\[
		P_\lambda(M)=M_\lambda.
		\]
		
		Because $T$ has pure point spectrum, its eigenspaces are mutually orthogonal and their closed linear span is $H$. Hence every $x\in H$ admits the orthogonal decomposition
		\[
		x=\sum_{\lambda\in\sigma_p(T)}P_\lambda x.
		\]
		For $x\in M$, each summand belongs to $P_\lambda(M)=M_\lambda\subset M$. Therefore
		\[
		M=
		\bigoplus_{\lambda\in\sigma_p(T)}M_\lambda.
		\]
	\end{proof}
	
	We can now prove the necessity of the reverse well-order condition on the eigenvalue moduli.
	
	\begin{lemma}\label{lem:rwo}
		Let $T$ be normal with pure point spectrum. If $T\in\theta(H)$, then
		\[
		A(T)=\{|\lambda|:\lambda\in\sigma_p(T)\}
		\]
		is reverse well-ordered.
	\end{lemma}
	
	\begin{proof}
		Suppose that $A(T)$ is not reverse well-ordered. By Definition~\ref{def:rwo} and the stated equivalence for subsets of $\mathbb{R}$, there exists a strictly increasing sequence of moduli; thus we can choose eigenvalues $\{\lambda_k\}$ such that
		\[
		|\lambda_1|<|\lambda_2|<\cdots.
		\]
		Set
		\[
		r=\lim_{k\to\infty}|\lambda_k|,
		\qquad
		M=
		\bigoplus_{k=1}^\infty
		\ker(T-\lambda_k I).
		\]
		Then $M\in\operatorname{Red}(T)$ and $T|_M$ is a diagonal normal operator with eigenvalue moduli $|\lambda_k|$. Therefore
		\[
		\omega(T|_M)
		=
		\|T|_M\|
		=
		\sup_k|\lambda_k|
		=
		r.
		\]
		Since $T\in\theta(H)$, there exists a unit vector $x\in M$ such that $|\langle Tx,x\rangle|=r$. Write
		\[
		x=\sum_{k=1}^\infty x_k,
		\qquad
		x_k\in\ker(T-\lambda_k I),
		\qquad
		\sum_{k=1}^\infty\|x_k\|^2=1.
		\]
		Then
		\[
		\langle Tx,x\rangle
		=
		\sum_{k=1}^\infty\lambda_k\|x_k\|^2,
		\]
		and hence
		\[
		\begin{aligned}
			r\leqslant
			\sum_{k=1}^\infty|\lambda_k|\|x_k\|^2=
			r-\sum_{k=1}^\infty(r-|\lambda_k|)\|x_k\|^2
			<r,
		\end{aligned}
		\]
		because $r-|\lambda_k|>0$ for all $k$ and at least one $\|x_k\|$ is nonzero. This contradiction proves the claim.
	\end{proof}
	
	Combining the two necessary conditions with the decomposition lemma, we arrive at the complete spectral characterization of normal operators in $\theta(H)$.
	
	\begin{theorem}\label{thm:normal}
		Let $H$ be a separable complex Hilbert space and let $T\in\mathcal{B}(H)$ be normal. Then
		\[
		T\in\theta(H)
		\]
		if and only if:
		\begin{enumerate}
			\item $T$ has pure point spectrum;
			\item $\{|\lambda|:\lambda\in\sigma_p(T)\}$ is reverse well-ordered.
		\end{enumerate}
	\end{theorem}
	
	\begin{proof}
		The necessity follows from Lemmas~\ref{lem:purepoint} and~\ref{lem:rwo}.
		
		Conversely, assume that $T$ has pure point spectrum and that
		\[
		A(T)=\{|\lambda|:\lambda\in\sigma_p(T)\}
		\]
		is reverse well-ordered. Let $0\neq M\in\operatorname{Red}(T)$. By Lemma~\ref{lem:reducing},
		\[
		M=
		\bigoplus_{\lambda\in\sigma_p(T)}M_\lambda,
		\qquad
		M_\lambda=M\cap\ker(T-\lambda I).
		\]
		Set
		\[
		S_M=
		\{|\lambda|:\lambda\in\sigma_p(T),\ M_\lambda\neq\{0\}\}.
		\]
		Since $M\neq\{0\}$ and $M=\overline{\bigoplus_\lambda M_\lambda}$, at least one $M_\lambda$ is nonzero; hence $S_M$ is nonempty. By the reverse well-ordering hypothesis, every nonempty subset of $A(T)$ has a maximum. Applying this to $S_M\subseteq A(T)$ gives a value $|\lambda_0|=\max S_M$. On each nonzero $M_\lambda$,
		\[
		T|_{M_\lambda}=\lambda I_{M_\lambda},
		\]
		\[
		\omega(T|_M)
		=
		\|T|_M\|
		=
		\max S_M
		=
		|\lambda_0|.
		\]
		Choose a unit vector $x_0\in M_{\lambda_0}$. Then
		\[
		Tx_0=\lambda_0x_0
		\]
		and hence
		\[
		|\langle Tx_0,x_0\rangle|
		=
		|\lambda_0|
		=
		\omega(T|_M).
		\]
		Since $M$ was arbitrary, $T\in\theta(H)$.
	\end{proof}
	
	The spectral characterization above depends essentially on normality: without normality, pure point spectrum is not necessary for membership in $\theta(H)$. The following example illustrates this phenomenon.
	
	\begin{example}\label{ex:nonnormal}
		Let $H=\ell^2$ and define the rank-one operator
		\[
		Te_2=e_1,\qquad Te_n=0\quad(n\neq2).
		\]
		Then $T$ is compact and non-normal. By Proposition~\ref{prop:compact}, $T\in\theta(H)$. However,
		\[
		\sigma_p(T)=\{0\},
		\qquad
		\ker T=\overline{\operatorname{span}}\{e_1,e_3,e_4,\ldots\}\neq H.
		\]
		Thus $T$ does not have pure point spectrum. Hence, for non-normal operators, membership in $\theta(H)$ does not imply pure point spectrum.
	\end{example}

	\section{Algebraic properties of \texorpdfstring{$\theta(H)$}{theta(H)}}\label{sec:algebra}

	Having obtained the spectral characterization for normal operators, we now turn to the basic algebraic properties of the class $\theta(H)$. The first four propositions below show that $\theta(H)$ behaves well under the most natural algebraic operations: nonzero scalar multiplication, taking adjoints, unitary similarity, and restriction to reducing subspaces.
	
	\begin{proposition}\label{prop:scalar}
		Let $T\in\mathcal{B}(H)$ and $\alpha\in\mathbb{C}\setminus\{0\}$. Then
		\[
		T\in\theta(H)
		\iff
		\alpha T\in\theta(H),
		\]
		and $0\in\theta(H)$.
	\end{proposition}
	
	\begin{proof}
		We have $\operatorname{Red}(\alpha T)=\operatorname{Red}(T)$, and for every nonzero reducing subspace $M$,
		\[
		\omega((\alpha T)|_M)
		=
		|\alpha|\,\omega(T|_M),
		\qquad
		|\langle\alpha Tx,x\rangle|
		=
		|\alpha|\,|\langle Tx,x\rangle|.
		\]
		Thus attainment is equivalent. The zero operator clearly belongs to $\theta(H)$.
	\end{proof}
	
	\begin{proposition}\label{prop:adjoint}
		For every $T\in\mathcal{B}(H)$,
		\[
		T\in\theta(H)
		\iff
		T^*\in\theta(H).
		\]
	\end{proposition}
	
	\begin{proof}
		We have $\operatorname{Red}(T)=\operatorname{Red}(T^*)$. If $0\neq M\in\operatorname{Red}(T)$, then
		\[
		(T|_M)^*=T^*|_M,
		\]
		and for every unit vector $x\in M$,
		\[
		|\langle T^*x,x\rangle|
		=
		|\overline{\langle Tx,x\rangle}|
		=
		|\langle Tx,x\rangle|.
		\]
		Hence
		\[
		\omega(T^*|_M)=\omega(T|_M),
		\]
		with identical attainment properties.
	\end{proof}
	
	\begin{proposition}\label{prop:unitary}
		Let $T\in\mathcal B(H)$ and let $U\in\mathcal B(H)$ be unitary. Then
		\[
		T\in\theta(H)
		\iff
		UTU^*\in\theta(H).
		\]
	\end{proposition}
	
	\begin{proof}
		We have
		\[
		M\in\operatorname{Red}(T)
		\iff
		UM\in\operatorname{Red}(UTU^*),
		\]
		and
		\[
		(UTU^*)|_{UM}=U(T|_M)U^*.
		\]
		Numerical radius is invariant under unitary similarity, and $U$ maps the unit sphere of $M$ onto the unit sphere of $UM$. Therefore the attainment property is preserved.
	\end{proof}
	
	\begin{proposition}\label{prop:restriction}
		Suppose that $T\in\theta(H)$ and $0\neq M\in\operatorname{Red}(T)$. Then
		\[
		T|_M\in\theta(M).
		\]
	\end{proposition}
	
	\begin{proof}
		Let $0\neq L\in\operatorname{Red}(T|_M)$. Since $M$ reduces $T$,
		\[
		(T|_M)^*=T^*|_M.
		\]
		Hence $L$ reduces $T$ and so $L\in\operatorname{Red}(T)$. Since $T\in\theta(H)$, there is a unit vector $x\in L$ such that
		\[
		|\langle Tx,x\rangle|
		=
		\omega(T|_L).
		\]
		As $L\subset M$,
		\[
		(T|_M)|_L=T|_L,
		\]
		and thus
		\[
		|\langle(T|_M)x,x\rangle|
		=
		\omega((T|_M)|_L).
		\]
		Therefore $T|_M\in\theta(M)$.
	\end{proof}
	
	Next we record a simple consequence of the fact that compact operators belong to $\theta(H)$: one-sided multiplication by a compact operator always produces an operator in $\theta(H)$.
	
	\begin{proposition}\label{prop:compact-mult}
		If $K\in\mathcal{B}(H)$ is compact and $T\in\mathcal{B}(H)$ is arbitrary, then
		\[
		TK,\ KT\in\theta(H).
		\]
	\end{proposition}
	
	\begin{proof}
		Both $TK$ and $KT$ are compact, so the assertion follows from Proposition~\ref{prop:compact}.
	\end{proof}
	
	The numerical radius satisfies the classical power inequality
	$\omega(T^n)\leqslant\omega(T)^n$ for every positive integer $n$
	\cite{pearcy1966}. To prove that positive powers preserve membership in
	$\theta(H)$ for normal operators, we use the spectral characterization in
	Theorem~\ref{thm:normal}. The same argument applies to the absolute value.
	Example~\ref{ex:nonnormal-power} will show that power stability can fail
	without normality.
	
	\begin{proposition}\label{prop:powers}
		Let $T$ be normal and assume $T\in\theta(H)$. Then, for every integer $n\geqslant1$,
		\[
		T^n\in\theta(H).
		\]
		Moreover,
		\[
		|T|=(T^*T)^{\frac12}\in\theta(H).
		\]
	\end{proposition}
	
	\begin{proof}
		By Theorem~\ref{thm:normal}, $T$ has pure point spectrum and
		\[
		A(T)=\{|\lambda|:\lambda\in\sigma_p(T)\}
		\]
		is reverse well-ordered. Hence $H$ has an orthonormal basis
		$\{e_j\}_{j\in J}$ consisting of eigenvectors of $T$, say
		$Te_j=\lambda_j e_j$. The same basis diagonalizes $T^n$, with
		$T^ne_j=\lambda_j^n e_j$. Consequently,
		\[
		\sigma_p(T^n)=\{\lambda_j^n:j\in J\}
		\quad\text{and}\quad
		A(T^n)=\{|\lambda|^n:\lambda\in\sigma_p(T)\}.
		\]
		Since $t\mapsto t^n$ is strictly increasing on $[0,\infty)$, the image of a reverse well-ordered set under this map is again reverse well-ordered. Theorem~\ref{thm:normal} therefore gives $T^n\in\theta(H)$.
		
		If $Te_j=\lambda_j e_j$, then
		\[
		|T|e_j=|\lambda_j|e_j.
		\]
		Thus $|T|$ is diagonal in the same orthonormal basis, has pure point
		spectrum, and satisfies
		\[
		A(|T|)=A(T).
		\]
		Hence $|T|\in\theta(H)$.
	\end{proof}
	
	The stability of $\theta(H)$ under finite orthogonal direct sums is proved
	in Section~\ref{sec:finite-direct-sum}, using graph-subspace methods.
	The required technical lemmas on closed unbounded operators are collected
	in Appendix~\ref{app:unbounded}.
	
	We conclude this section with an example showing that a natural operation does not preserve membership in $\theta(H)$: translation by a scalar multiple of the identity may destroy the reverse well-order of the eigenvalue moduli and hence membership. This example also shows that $\theta(H)$ is neither a linear space nor translation-invariant. (The failure of infinite orthogonal direct sums is postponed to Section~\ref{sec:finite-direct-sum}, after the proof of the finite direct-sum theorem.)
	
	\begin{example}\label{ex:translation}
		Let $H=\ell^2(\mathbb N)$ with standard orthonormal basis $(e_n)$, and define
		\[
		Te_n=\left(1+\frac1n\right)e_n.
		\]
		Then $T$ is normal and $T\in\theta(H)$. On the other hand, $-2I\in\theta(H)$, but
		\[
		(T-2I)e_n=-\left(1-\frac1n\right)e_n.
		\]
		The moduli of the eigenvalues are strictly increasing and converge to $1$. Hence, by Theorem~\ref{thm:normal},
		\[
		T-2I\notin\theta(H).
		\]
		Indeed,
		\[
		\langle(T-2I)x,x\rangle
		=
		-1+\sum_n\frac1n|x_n|^2>-1
		\]
		for every unit vector $x=\sum_nx_ne_n$. Since $\frac1n\leqslant1$, the displayed expression lies in $(-1,0]$, so
		\[
		|\langle(T-2I)x,x\rangle|<1
		\]
		for every unit vector $x$. On the other hand, $T-2I$ is self-adjoint with spectrum $\left\{-1+\frac1n:n\geqslant1\right\}\cup\{-1\}$, and therefore
		\[
		\omega(T-2I)=\sup\{|\lambda|:\lambda\in\sigma(T-2I)\}=1.
		\]
		Thus the supremum is not attained. This shows directly that $T-2I\notin\theta(H)$, and hence that $\theta(H)$ is generally neither a linear space nor translation-invariant.
	\end{example}

	\section{Finite orthogonal direct sums}\label{sec:finite-direct-sum}

	Although the class $\theta(H)$ is defined for bounded operators, its stability
	under finite orthogonal direct sums is most conveniently proved by decomposing a
	reducing subspace into the graph of a closed, possibly unbounded, operator.
	We begin with an explicit non-normal example that illustrates the proof
	without the technical difficulties of unbounded domains. We then outline
	the proof and introduce the graph-subspace framework before proving the
	general theorem. The technical results on closed operators used in the
	argument are proved in Appendix~\ref{app:unbounded}, following
	\cite[Chapter X]{Conway1990} and \cite{kato1976}.
	Remark~\ref{rem:unbounded-graph} explains why unbounded graph operators
	are necessary, and Example~\ref{ex:infinite-sum} gives the infinite-sum
	obstruction. Throughout this section, $H_1,H_2$ are separable complex
	Hilbert spaces.
	
	\subsection{A worked example: a non-normal finite direct sum}
	
	\begin{example}\label{ex:finite-direct-sum-worked}
		Let $\mathcal H=\ell^2(\mathbb N)$, put $H_0=\mathcal H\oplus\mathcal H$,
		and define $R\in\mathcal B(H_0)$ by
		\[
		R(x,y)=(y,0)\qquad(x,y\in\mathcal H).
		\]
		The operator $R$ is neither normal nor compact: $R^2=0$,
		$RR^*\neq R^*R$, and its restriction from $0\oplus\mathcal H$
		onto $\mathcal H\oplus0$ is an isometry. We show directly that
		\[
		R\in\theta(H_0),\qquad
		T:=R\oplus R\in\theta(H_0\oplus H_0).
		\]
		
		First consider a closed subspace $D\subset\mathcal H$. The subspace
		$D\oplus D$ reduces $R$, and the restriction is given by
		$R_D(x,y)=(y,0)$ for $x,y\in D$. For a unit vector $(x,y)$,
		\[
		|\langle R_D(x,y),(x,y)\rangle|
		=|\langle y,x\rangle|
		\leqslant\|x\|\,\|y\|
		\leqslant\frac{\|x\|^2+\|y\|^2}{2}
		=\frac12.
		\]
		If $D\neq\{0\}$, equality holds at $\frac{(u,u)}{\sqrt{2}}$ for every unit vector $u\in D$. Thus $\omega(R_D)=\frac12$, and it is attained.
		In fact, these are \emph{all} the reducing subspaces of $R$:
		if an orthogonal projection $P$ on $H_0$ is written as
		\[
		P=\begin{pmatrix}P_{11}&P_{12}\\P_{21}&P_{22}\end{pmatrix},
		\]
		then $PR=RP$ and $PR^*=R^*P$ imply
		\[
		P_{12}=P_{21}=0,\qquad P_{11}=P_{22}=:Q.
		\]
		Since $P$ is an orthogonal projection, so is $Q$, and hence
		$\operatorname{ran}P=D\oplus D$ with $D=\operatorname{ran}Q$.
		Consequently $R\in\theta(H_0)$.
		
		Before examining arbitrary reducing subspaces of $T$, we illustrate
		the three pieces that occur in the general proof. Choose an
		orthogonal decomposition
		$\mathcal H=F\oplus E\oplus G$ with $F,E,G$ nonzero closed subspaces,
		and set
		\[
		N=F\oplus F,\quad K_1=E\oplus E,\quad
		L=E\oplus E,\quad L_0=G\oplus G,
		\]
		all regarded as subspaces of $H_0$. Define
		$S:K_1\to L$ by $Sx=2x$ and put
		\[
		M_g=\{(x,Sx):x\in K_1\},\qquad
		M=(0\oplus L_0)\oplus(N\oplus0)\oplus M_g.
		\]
		The three terms are mutually orthogonal. Both $R$ and $R^*$
		preserve $N$, $K_1=L$, and $L_0$, while
		\[
		SR|_{K_1}=R|_L S,\qquad
		SR^*|_{K_1}=R^*|_L S.
		\]
		It follows that each term reduces $T=R\oplus R$, so $M$ is a
		nonzero closed reducing subspace of $T$. The graph component is
		genuinely coupled across the two coordinates: if $0\neq x\in K_1$,
		then $(x,2x)\in M$ but $(x,0)\notin M$.
		
		We now follow the mechanism used later for an arbitrary graph
		component. Here $S^*S=4I_{K_1}$, its polar unitary is the identity
		from $K_1$ to $L$, and its graph inner product is
		\[
		\langle x,y\rangle_S
		=\langle x,y\rangle+\langle Sx,Sy\rangle
		=5\langle x,y\rangle.
		\]
		Thus
		\[
		V_g:K_1\longrightarrow M_g,\qquad
		V_gx=\frac{(x,2x)}{\sqrt5}
		\]
		is unitary and satisfies
		\[
		T|_{M_g}V_g=V_g(R|_{K_1}).
		\]
		Since $K_1=E\oplus E$, our calculation for $R_D$ yields
		\[
		\omega(T|_{M_g})=\omega(R|_{K_1})=\frac12,
		\]
		attained at $V_g\frac{(u,u)}{\sqrt{2}}$ for any unit vector $u\in E$.
		The restrictions to $0\oplus L_0$ and $N\oplus0$ have the same
		numerical radius $\frac12$, attained as well. Therefore
		\[
		\omega(T|_M)
		=\max\left\{\omega(R|_{L_0}),\omega(R|_N),
		\omega(T|_{M_g})\right\}
		=\frac12,
		\]
		and this numerical radius is attained. This is the pattern of the
		general proof: separate the vertical, horizontal and graph parts;
		transfer attainment through the graph norm; then take the maximum
		of finitely many attained numerical radii.
		
		Finally, the conclusion $T\in\theta(H_0\oplus H_0)$ holds for
		\emph{every} reducing subspace, not just the one exhibited above.
		Indeed, the coordinate-reordering unitary
		\[
		\mathcal U((x_1,y_1),(x_2,y_2))
		=((x_1,x_2),(y_1,y_2))
		\]
		identifies $T$ with
		\[
		\mathcal U T\mathcal U^*
		=\begin{pmatrix}0&I_{\mathcal H\oplus\mathcal H}\\
			0&0\end{pmatrix}.
		\]
		The same projection calculation as above shows that an arbitrary
		nonzero reducing subspace of $T$ has the form
		$\mathcal U^*(D\oplus D)$ for some nonzero closed
		$D\subset\mathcal H\oplus\mathcal H$. Its restriction again has
		numerical radius $\frac12$, attained at
		$\mathcal U^*\frac{(v,v)}{\sqrt{2}}$ for a unit vector $v\in D$.
		Hence $R\oplus R\in\theta(H_0\oplus H_0)$, as claimed.
	\end{example}
	
	\subsection{Graph-subspace framework and proof roadmap}\label{subsec:graph-roadmap}
	
	A linear operator $S:D(S)\subset H_1\to H_2$ is \emph{closed} if its graph
	\[
	\operatorname{graph}(S)=\{(x,Sx):x\in D(S)\}
	\]
	is closed in $H_1\oplus H_2$. The \emph{graph inner product} and \emph{graph norm}
	on $D(S)$ are
	\[
	\langle x,y\rangle_S=\langle x,y\rangle+\langle Sx,Sy\rangle,
	\qquad
	\|x\|_S^2=\|x\|^2+\|Sx\|^2.
	\]
	The adjoint and positive-square-root facts, together with self-contained
	proofs of the graph-completion, polar-intertwining, and graph-norm
	extension lemmas, are collected in Appendix~\ref{app:unbounded}.
	Specifically, Lemma~\ref{lem:graph-completion} identifies the graph with
	the graph-norm Hilbert space, Lemma~\ref{lem:polar-intertwining} gives
	the unitary equivalence of the coordinate restrictions, and
	Lemma~\ref{lem:graph-bounded-extension} transfers their self-adjoint
	parts to the graph-norm space.
	
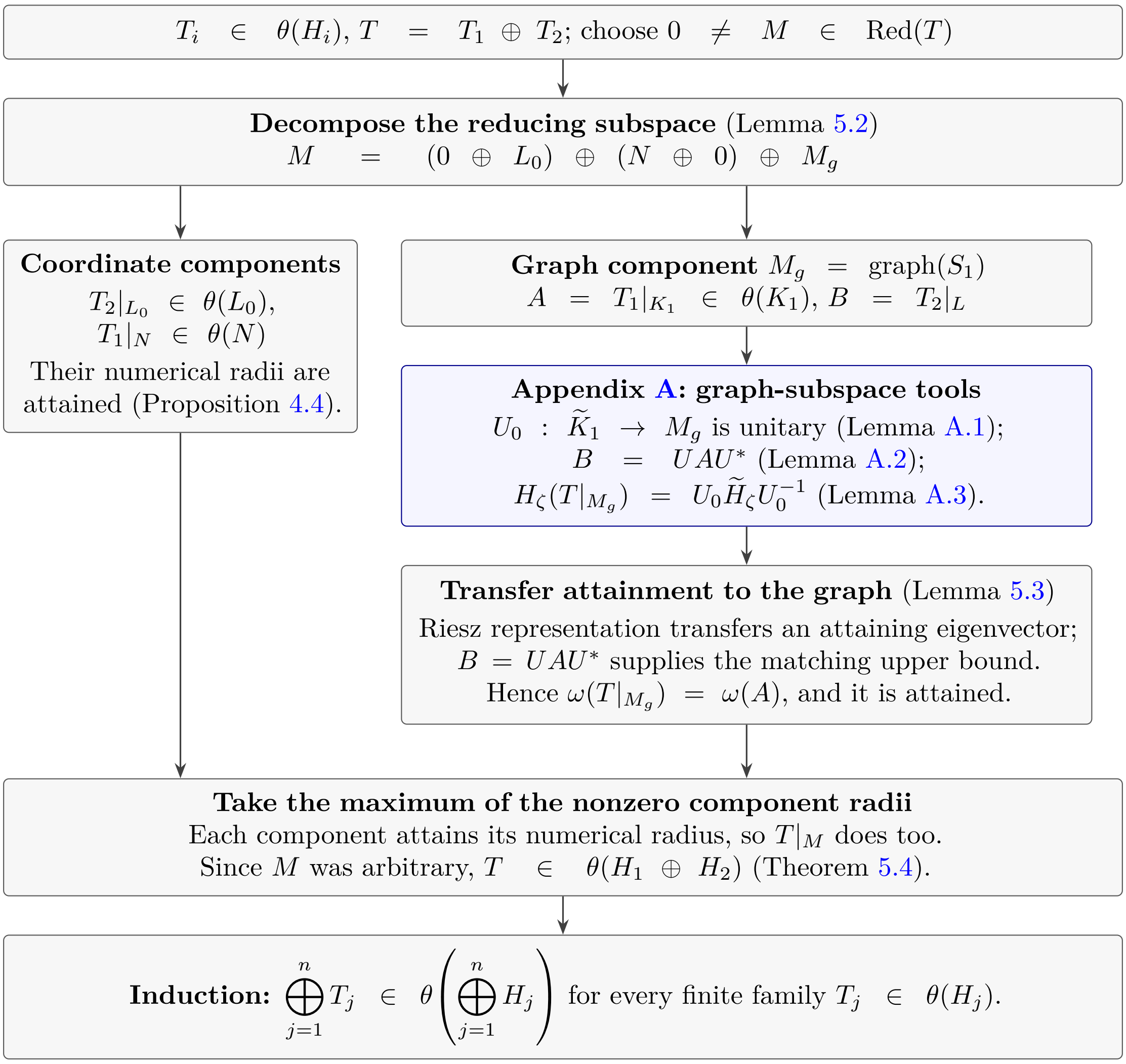
\begin{figure}[H]
		\centering
		\begin{tikzpicture}[
			font=\small,
			roadmapbox/.style={draw=black!65,rounded corners=2pt,align=center,
				inner xsep=6pt,inner ysep=5pt,text width=14.2cm},
			roadmapstep/.style={roadmapbox,fill=black!3},
			roadmaptool/.style={roadmapbox,fill=blue!4,draw=blue!55!black},
			roadmapflow/.style={-{Stealth[length=2mm]},semithick,draw=black!75}
			]
			\node[roadmapstep] (start) {$T_i\in\theta(H_i)$, $T=T_1\oplus T_2$; choose
				$0\ne M\in\operatorname{Red}(T)$};
			\node[roadmapstep,below=5mm of start] (split) {
				\textbf{Decompose the reducing subspace}
				(Lemma~\ref{lem:direct-sum-decomposition})\\
				$M=(0\oplus L_0)\oplus(N\oplus0)\oplus M_g$};
			\node[roadmapstep,text width=4.2cm,anchor=north] (coordinate)
			at ([xshift=-5cm,yshift=-7mm]split.south) {
				\textbf{Coordinate components}\\[2pt]
				$T_2|_{L_0}\in\theta(L_0)$,\\
				$T_1|_N\in\theta(N)$\\[2pt]
				Their numerical radii are attained
				(Proposition~\ref{prop:restriction}).};
			\node[roadmapstep,text width=8.6cm,anchor=north] (graph)
			at ([xshift=2.4cm,yshift=-7mm]split.south) {
				\textbf{Graph component} $M_g=\operatorname{graph}(S_1)$\\
				$A=T_1|_{K_1}\in\theta(K_1)$, $B=T_2|_L$};
			\node[roadmaptool,text width=8.6cm,below=5mm of graph] (tools) {
				\textbf{Appendix~\ref{app:unbounded}: graph-subspace tools}\\[2pt]
				$U_0:\widetilde K_1\to M_g$ is unitary
				(Lemma~\ref{lem:graph-completion});\\
				$B=UAU^*$ (Lemma~\ref{lem:polar-intertwining});\\
				$H_\zeta(T|_{M_g})=U_0\widetilde H_\zeta U_0^{-1}$
				(Lemma~\ref{lem:graph-bounded-extension}).};
			\node[roadmapstep,text width=8.6cm,below=5mm of tools] (attainment) {
				\textbf{Transfer attainment to the graph}
				(Lemma~\ref{lem:graph-attainment})\\[2pt]
				Riesz representation transfers an attaining eigenvector;\\
				$B=UAU^*$ supplies the matching upper bound.\\
				Hence $\omega(T|_{M_g})=\omega(A)$, and it is attained.};
			\node[roadmapstep,anchor=north] (maximum)
			at ([xshift=-2.4cm,yshift=-7mm]attainment.south) {
				\textbf{Take the maximum of the nonzero component radii}\\
				Each component attains its numerical radius, so $T|_M$ does too.\\
				Since $M$ was arbitrary, $T\in\theta(H_1\oplus H_2)$
				(Theorem~\ref{thm:finite-direct-sum}).};
			\node[roadmapstep,below=5mm of maximum] (induction) {
				\textbf{Induction:}
				$\displaystyle\bigoplus_{j=1}^{n}T_j\in
				\theta\!\left(\bigoplus_{j=1}^{n}H_j\right)$
				for every finite family $T_j\in\theta(H_j)$.};
			\draw[roadmapflow] (start) -- (split);
			\draw[roadmapflow] ([xshift=-5cm]split.south) -- (coordinate.north);
			\draw[roadmapflow] ([xshift=2.4cm]split.south) -- (graph.north);
			\draw[roadmapflow] (graph) -- (tools);
			\draw[roadmapflow] (tools) -- (attainment);
			\draw[roadmapflow] (coordinate.south) -- (coordinate.south |- maximum.north);
			\draw[roadmapflow] (attainment.south) -- (attainment.south |- maximum.north);
			\draw[roadmapflow] (maximum) -- (induction);
		\end{tikzpicture}
		\caption{Proof roadmap for stability under finite orthogonal direct sums.
			The appendix supplies the technical tools for the graph branch.}
		\label{fig:finite-direct-sum-roadmap}
	\end{figure}

	It suffices to prove the theorem for two summands and then use induction.
	Figure~\ref{fig:finite-direct-sum-roadmap} outlines the two-summand argument
	carried out in Section~\ref{subsec:finite-direct-sum-proof}. The key step is
	to treat the graph component of an arbitrary reducing subspace: the
	appendix lemmas allow an attaining eigenvector to be transferred by Riesz
	representation in the graph norm. The coordinate components are handled
	by the restriction property of $\theta$, and a finite maximum combines
	the components. Zero components are omitted throughout the argument.
	
	For reference, Table~\ref{tab:notation} records the spaces and operators
	used in the graph argument and the finite direct-sum decomposition.
	
	\begin{longtable}{@{}>{$}l<{$}p{0.7\textwidth}@{}}
		\caption{Notation table for the graph-subspace decomposition and unbounded operator arguments.}\label{tab:notation}\\
		\toprule
		\text{\textbf{Symbol}} & \textbf{Meaning / Space} \\
		\midrule
		\endfirsthead
		\multicolumn{2}{c}{{\tablename\ \thetable{} -- continued from previous page}} \\
		\toprule
		\text{\textbf{Symbol}} & \textbf{Meaning / Space} \\
		\midrule
		\endhead
		\midrule
		\multicolumn{2}{r}{{Continued on next page}} \\
		\endfoot
		\bottomrule
		\endlastfoot
		H_1, H_2 & The complex Hilbert spaces for the direct sum $H_1 \oplus H_2$. \\
		T_1, T_2 & Bounded operators in $\mathcal{B}(H_1)$ and $\mathcal{B}(H_2)$. \\
		M & A nonzero closed reducing subspace of $T_1 \oplus T_2$. \\
		L_0 & A closed reducing subspace of $T_2$ such that $0 \oplus L_0 = M \cap (0 \oplus H_2)$. \\
		N & The kernel of $S_0$; a closed reducing subspace of $T_1$ in $H_1$. \\
		S_0 & A closed, single-valued operator from $D(S_0) \subset H_1$ to $H_2$ whose graph is $M \ominus (0 \oplus L_0)$. \\
		K & The closure $\overline{D(S_0)}$ in $H_1$. \\
		L & The closure of the range $\overline{\operatorname{ran}S_0}$ in $H_2$. \\
		K_1 & The orthogonal complement $K \ominus N$, on which the operator is injective. \\
		S_1 & The restriction $S_0|_{D(S_1)}$, which is a closed injective operator with dense domain in $K_1$ and dense range in $L$. \\
		D(S_1) & The domain of $S_1$, defined as $D(S_0) \cap K_1$. \\
		M_g & The graph subspace $\{(x, S_1x) : x \in D(S_1)\} \subset H_1 \oplus H_2$. \\
		W & The positive self-adjoint operator $S_1^*S_1$ on $K_1$. \\
		\widetilde{K}_1 & The completion of $D(S_1)$ with respect to the graph norm $\|\cdot\|_{S_1}$. \\
		U_0 & The unitary operator mapping $\widetilde{K}_1$ onto $M_g$, given by $U_0x = (x, S_1x)$. \\
		U & The unitary operator in the polar decomposition $S_1 = UW^{\frac12}$, mapping $K_1$ to $L$. \\
	\end{longtable}
	
	\subsection{Finite orthogonal direct sums}\label{subsec:finite-direct-sum-proof}
	
	We now carry out the proof outlined in
	Figure~\ref{fig:finite-direct-sum-roadmap}, using the technical lemmas from
	Appendix~\ref{app:unbounded}. The first step describes the structure of
	an arbitrary reducing subspace of $T_1\oplus T_2$.
	
	\begin{lemma}[Decomposition of reducing subspaces]\label{lem:direct-sum-decomposition}
		Let $T=T_1\oplus T_2$ on $H_1\oplus H_2$, where
		$T_i\in\mathcal B(H_i)$, and let $0\neq M\in\operatorname{Red}(T)$. Then $M$
		admits an orthogonal decomposition
		\[
		M=(0\oplus L_0)\oplus(N\oplus0)\oplus M_g,
		\]
		where
		\begin{enumerate}
			\item $L_0\in\operatorname{Red}(T_2)$ and $N\in\operatorname{Red}(T_1)$;
			\item $M_g=\{(x,S_1x):x\in D(S_1)\}$ is the graph of a closed injective
			operator $S_1:D(S_1)\subset K_1\to L$ with dense domain in a closed
			subspace $K_1\in\operatorname{Red}(T_1)$ and dense range in a closed
			subspace $L\in\operatorname{Red}(T_2)$;
			\item Writing $A=T_1|_{K_1}$ and $B=T_2|_L$, one has
			$S_1Ax=BS_1x$ and $S_1A^*x=B^*S_1x$ for every $x\in D(S_1)$;
			\item $B=UAU^*$ for a unitary operator $U:K_1\to L$.
		\end{enumerate}
	\end{lemma}
	
	\begin{proof}
		Put
		\[
		M_0=M\cap(0\oplus H_2).
		\]
		Since $M$ reduces $T$, we have $M_0=0\oplus L_0$ with
		$L_0\in\operatorname{Red}(T_2)$. Set $M_1=M\ominus M_0$. Then $M_0$ and
		$M_1$ both reduce $T$, and
		\[
		M_1\cap(0\oplus H_2)=\{0\}.
		\]
		Hence the first-coordinate projection $P_1$ is injective on $M_1$. Let
		\[
		D(S_0)=P_1(M_1),
		\]
		and for $x\in D(S_0)$ define $S_0x=y$, where $y\in H_2$ is the unique vector
		such that $(x,y)\in M_1$. Then $S_0$ is a closed single-valued operator and
		\[
		M_1=\{(x,S_0x):x\in D(S_0)\}.
		\]
		Put $K=\overline{D(S_0)}$ and $L=\overline{\operatorname{ran}S_0}$.
		Orthogonality to $0\oplus L_0$ gives $L\subset L_0^\perp$. Since
		$M_1$ reduces $T$, for $x\in D(S_0)$ we have
		\[
		T_1x\in D(S_0),\qquad S_0T_1x=T_2S_0x,
		\]
		and
		\[
		T_1^*x\in D(S_0),\qquad S_0T_1^*x=T_2^*S_0x.
		\]
		Thus $K\in\operatorname{Red}(T_1)$ and $L\in\operatorname{Red}(T_2)$.
		
		Let $N=\ker S_0$. Since $S_0$ is closed, $N$ is closed; the intertwining
		relations give $T_1N\subset N$ and $T_1^*N\subset N$, so
		$N\in\operatorname{Red}(T_1)$. Set
		\[
		K_1=K\ominus N,\qquad
		D(S_1)=D(S_0)\cap K_1,\qquad
		S_1=S_0|_{D(S_1)}.
		\]
		Let $P_N$ be the orthogonal projection of $K$ onto $N$. Since
		$N\subset D(S_0)$, for every $x\in D(S_0)$ we have
		$P_Nx\in N\subset D(S_0)$ and
		$(I-P_N)x\in D(S_0)\cap K_1=D(S_1)$. Hence
		\[
		D(S_0)=N\oplus D(S_1),
		\qquad
		D(S_1)=(I-P_N)D(S_0).
		\]
		As $D(S_0)$ is dense in $K$ and $I-P_N$ maps $K$ continuously onto
		$K_1$, the space $D(S_1)$ is dense in $K_1$. The operator $S_1$ is
		closed because its graph is the intersection of the closed graph of
		$S_0$ with $K_1\oplus H_2$, and it is injective by construction.
		Moreover, $S_0$ vanishes on $N$, so
		$\operatorname{ran}S_1=\operatorname{ran}S_0$; in particular,
		$\operatorname{ran}S_1$ is dense in $L$. Since $K$ and $N$ reduce
		$T_1$, their orthogonal difference $K_1$ also reduces $T_1$.
		Writing $A=T_1|_{K_1}$ and $B=T_2|_L$, the intertwining relations
		become
		\[
		S_1Ax=BS_1x,\qquad
		S_1A^*x=B^*S_1x
		\qquad(x\in D(S_1)).
		\]
		Since $S_0$ vanishes on $N$ and $D(S_0)=N\oplus D(S_1)$, the graph of
		$S_0$ decomposes as
		\[
		M_1=(N\oplus0)\oplus M_g.
		\]
		Consequently
		\[
		M=(0\oplus L_0)\oplus(N\oplus0)\oplus M_g,
		\]
		where
		\[
		M_g=\{(x,S_1x):x\in D(S_1)\},
		\]
		and the three summands are mutually orthogonal and reduce $T$.
		
		It remains to prove (4). Put $W=S_1^*S_1$, and write
		$A=T_1|_{K_1}$ and $B=T_2|_L$. By
		Lemma~\ref{lem:polar-intertwining}, $S_1=UW^{\frac12}$ for a unitary
		$U:K_1\to L$, and the intertwining relations give $UA=BU$.
		Therefore $B=UAU^*$, as required.
	\end{proof}
	
	\begin{lemma}[Attainment on a graph subspace]\label{lem:graph-attainment}
		Under the notation of Lemma~\ref{lem:direct-sum-decomposition}, assume that
		$K_1\neq\{0\}$ and $T_1|_{K_1}\in\theta(K_1)$. Then
		$(T_1\oplus T_2)|_{M_g}$ attains its numerical radius and
		\[
		\omega\bigl((T_1\oplus T_2)|_{M_g}\bigr)
		=
		\omega(T_1|_{K_1}).
		\]
	\end{lemma}
	
	\begin{proof}
		Let $A=T_1|_{K_1}$ and $B=T_2|_L$. Let
		$U_0:\widetilde K_1\to M_g$ be the unitary from
		Lemma~\ref{lem:graph-completion}. By
		Lemma~\ref{lem:graph-bounded-extension}, for every $\zeta\in\mathbb{T}$,
		\[
		H_\zeta\bigl((T_1\oplus T_2)|_{M_g}\bigr)
		=U_0\widetilde H_\zeta U_0^{-1},
		\]
		where $\widetilde H_\zeta$ is the graph-norm extension associated
		with $A$ and $B$.
		Therefore
		\[
		\omega\bigl((T_1\oplus T_2)|_{M_g}\bigr)
		=
		\max_{\zeta\in\mathbb{T}}\sup\sigma(\widetilde H_\zeta).
		\]
		
		Since $A=T_1|_{K_1}\in\theta(K_1)$, there exist $\zeta_0\in\mathbb{T}$
		and a unit vector $x_0\in K_1$ such that
		\[
		H_{\zeta_0}(A)x_0=\omega(A)x_0.
		\]
		For $x\in D(S_1)$,
		\[
		|\langle x,x_0\rangle|\leqslant\|x\|\leqslant\|x\|_{S_1},
		\]
		so the functional $x\mapsto\langle x,x_0\rangle$ is continuous with respect to
		the graph norm and extends to a continuous linear functional on
		$\widetilde K_1$. By the Riesz representation theorem, there exists a unique
		$\tilde x_0\in\widetilde K_1$ such that
		\[
		\langle x,x_0\rangle=\langle x,\tilde x_0\rangle_{S_1}
		\quad(x\in D(S_1)).
		\]
		Since $D(S_1)$ is dense in $K_1$ and $x_0\neq0$, the functional $x\mapsto\langle x,x_0\rangle$ is not identically zero on $D(S_1)$. Hence its Riesz representative satisfies $\tilde x_0\neq0$. For $x\in D(S_1)$,
		\[
		\langle x,\widetilde H_{\zeta_0}\tilde x_0\rangle_{S_1}
		=\langle H_{\zeta_0}(A)x,x_0\rangle
		=\omega(A)\langle x,\tilde x_0\rangle_{S_1}.
		\]
		By density,
		\[
		\widetilde H_{\zeta_0}\tilde x_0
		=\omega(A)\tilde x_0.
		\]
		Hence
		\[
		\sup\sigma(\widetilde H_{\zeta_0})\geqslant\omega(A).
		\]
		
		Conversely, for $\zeta\in\mathbb{T}$ and $x\in D(S_1)$,
		\[
		\langle\widetilde H_\zeta x,x\rangle_{S_1}
		=\langle H_\zeta(A)x,x\rangle
		+\langle H_\zeta(B)S_1x,S_1x\rangle
		\leqslant \alpha_\zeta\|x\|_{S_1}^2,
		\]
		where
		\[
		\alpha_\zeta=
		\max\left\{
		\sup\sigma(H_\zeta(A)),
		\sup\sigma(H_\zeta(B))
		\right\}.
		\]
		By Lemma~\ref{lem:polar-intertwining}, $B=UAU^*$, so the two
		suprema coincide; denote their common value by $\alpha_\zeta$. Thus
		\[
		\sup\sigma(\widetilde H_\zeta)\leqslant\alpha_\zeta
		\leqslant\omega(A).
		\]
		Taking the maximum over $\zeta$ gives
		\[
		\omega\bigl((T_1\oplus T_2)|_{M_g}\bigr)
		\leqslant\omega(A).
		\]
		Combining the two inequalities yields
		\[
		\omega\bigl((T_1\oplus T_2)|_{M_g}\bigr)
		=\omega(A).
		\]
		To obtain an attaining unit vector, set
		\[
		y_0=U_0\left(\frac{\tilde x_0}{\|\tilde x_0\|_{S_1}}\right)\in M_g.
		\]
		Since $U_0$ is unitary, $\|y_0\|=1$, and the eigenvector relation gives
		\[
		H_{\zeta_0}\bigl((T_1\oplus T_2)|_{M_g}\bigr)y_0=\omega(A)y_0.
		\]
		Consequently,
		\[
		\begin{aligned}
			\omega(A)
			&=\operatorname{Re}\bigl(\zeta_0\langle(T_1\oplus T_2)y_0,y_0\rangle\bigr)\\
			&\leqslant |\langle(T_1\oplus T_2)y_0,y_0\rangle|
			\leqslant\omega\bigl((T_1\oplus T_2)|_{M_g}\bigr)=\omega(A).
		\end{aligned}
		\]
		Thus $y_0$ attains the numerical radius on $M_g$.
	\end{proof}
	
	\begin{theorem}[Finite direct sums]\label{thm:finite-direct-sum}
		Let $H_1,H_2$ be separable complex Hilbert spaces and let $T_i\in\mathcal B(H_i)$
		satisfy
		\[
		T_i\in\theta(H_i)\qquad(i=1,2).
		\]
		Then
		\[
		T_1\oplus T_2\in\theta(H_1\oplus H_2).
		\]
	\end{theorem}
	
	\begin{proof}
		Let $T=T_1\oplus T_2$ and take $0\neq M\in\operatorname{Red}(T)$. By
		Lemma~\ref{lem:direct-sum-decomposition},
		\[
		M=(0\oplus L_0)\oplus(N\oplus0)\oplus M_g.
		\]
		The summands are mutually orthogonal and reduce $T$.
		
		If $L_0\neq\{0\}$, then Proposition~\ref{prop:restriction} gives
		$T_2|_{L_0}\in\theta(L_0)$, so $T|_{0\oplus L_0}=T_2|_{L_0}$ attains its
		numerical radius; if $L_0=\{0\}$, this summand is absent. Similarly, if
		$N\neq\{0\}$, then $T|_{N\oplus0}=T_1|_N$ attains its numerical radius
		because $T_1|_N\in\theta(N)$; if $N=\{0\}$, this summand is absent.
		
		For the graph summand, if $M_g=\{0\}$ there is nothing to prove. Otherwise
		$K_1\neq\{0\}$, and Proposition~\ref{prop:restriction} gives
		$T_1|_{K_1}\in\theta(K_1)$. By Lemma~\ref{lem:graph-attainment},
		$(T_1\oplus T_2)|_{M_g}$ attains its numerical radius.
		
		Thus $T|_M$ is the orthogonal direct sum of at most three nonzero
		restrictions, each of which attains its numerical radius. For a finite
		orthogonal direct sum $A_1\oplus\cdots\oplus A_m$, one has
		\[
		\omega(A_1\oplus\cdots\oplus A_m)
		=\max_{1\leqslant j\leqslant m}\omega(A_j),
		\]
		because, for $x=x_1\oplus\cdots\oplus x_m$ with $\|x\|=1$,
		\[
		|\langle(A_1\oplus\cdots\oplus A_m)x,x\rangle|
		\leqslant\sum_{j=1}^m\omega(A_j)\|x_j\|^2
		\leqslant\max_j\omega(A_j),
		\]
		while the reverse inequality follows by testing unit vectors in each
		individual summand. Hence a unit vector attaining the largest component
		radius gives attainment on $M$. Since $M\in\operatorname{Red}(T)$ was arbitrary,
		\[
		T_1\oplus T_2\in\theta(H_1\oplus H_2).
		\]
	\end{proof}
	
	\begin{remark}[Why unbounded graph operators are necessary]\label{rem:unbounded-graph}
		The graph operators in the proof of Theorem~\ref{thm:finite-direct-sum}
		cannot, in general, be assumed bounded. This issue occurs even
		for two identical positive compact summands. Let
		$H_0=\ell^2(\mathbb N)$ with its standard orthonormal basis
		$\{e_n\}_{n\geqslant1}$, and define
		\[
		Ae_n=\frac1n e_n.
		\]
		Since $A$ is compact, Proposition~\ref{prop:compact} gives
		$A\in\theta(H_0)$. Consider the densely defined operator
		\[
		Se_n=ne_n,\qquad
		D(S)=\left\{x=(x_n)\in\ell^2(\mathbb N):
		\sum_{n=1}^{\infty}n^2|x_n|^2<\infty\right\}.
		\]
		The operator $S$ is positive, self-adjoint, closed and unbounded.
		Moreover, $\operatorname{ran}A=D(S)$ and
		\[
		SAx=x\quad(x\in H_0),\qquad
		ASx=x\quad(x\in D(S)).
		\]
		Hence for every $x\in D(S)$,
		\[
		(A\oplus A)(x,Sx)
		=(Ax,ASx)=(Ax,x)=(Ax,SAx).
		\]
		Its graph
		\[
		M=\operatorname{graph}(S)
		=\{(x,Sx):x\in D(S)\}
		\]
		is therefore invariant under $A\oplus A$; since $A\oplus A$
		is self-adjoint, $M$ reduces $A\oplus A$.
		
		The first-coordinate projection of $M$ is precisely $D(S)$,
		which is dense but not closed in $H_0$. Indeed,
		$x=\left\{\frac1n\right\}_{n\geqslant1}\in\ell^2$ does not belong to $D(S)$,
		although its finite truncations belong to $D(S)$ and converge
		to $x$. Thus $M$ is not a coordinatewise orthogonal direct sum
		$M_1\oplus M_2$ of closed subspaces of $H_0$, and its graph
		cannot be represented by a bounded operator on a closed
		first-coordinate domain. Nevertheless,
		\[
		\omega\bigl((A\oplus A)|_M\bigr)=1
		\]
		is attained at $\frac{(e_1,e_1)}{\sqrt{2}}$. The example shows precisely
		why the decomposition and graph-attainment arguments above must
		allow \emph{closed unbounded} operators: bounded graph maps would
		not cover all reducing subspaces, even for summands in $\theta$.
	\end{remark}

	By induction, any finite orthogonal direct sum $T_1\oplus\cdots\oplus T_n$ with
	$T_i\in\theta(H_i)$ belongs to $\theta(H_1\oplus\cdots\oplus H_n)$.
	
	\begin{remark}\label{rem:many-reducing}
		The finite direct-sum theorem may come as a surprise when compared with the
		rich structure of reducing subspaces of direct sums. Indeed, it is known from
		\cite{FangWuJiang2003} that the number of reducing subspaces of an operator
		is either finite or uncountably infinite, and that an operator has only
		finitely many reducing subspaces exactly when it is the direct sum of
		finitely many pairwise unitarily inequivalent irreducible operators. In
		particular, if $A$ and $B$ are irreducible and unitarily equivalent, then
		$A\oplus B$ has uncountably many reducing subspaces; more generally,
		$A^{(n)}$ has a reducing subspace lattice isomorphic to that of $I_n$. Thus a
		finite orthogonal direct sum can carry a very large and complicated family of reducing
		subspaces. Nevertheless, Theorem~\ref{thm:finite-direct-sum} shows that
		$\theta(H)$ is inherited by finite orthogonal direct sums without any
		normality assumption. It is striking that, despite this
		abundance of reducing subspaces, the attainment of the numerical radius on
		every one of them is preserved under finite orthogonal direct sums.
	\end{remark}
	
	The finite direct-sum theorem stands in contrast with the infinite case. The
	following example shows that infinite orthogonal direct sums need not preserve
	membership in $\theta(H)$.
	
	\begin{example}\label{ex:infinite-sum}
		Let
		\[
		T_n=\left(1-\frac1n\right)I_{\mathbb{C}}.
		\]
		Then $T_n\in\theta(\mathbb{C})$ for every $n\in\mathbb{N}$. However, the infinite orthogonal direct sum
		\[
		T=\bigoplus_{n=1}^\infty T_n
		\]
		satisfies
		\[
		Te_n=\left(1-\frac1n\right)e_n,
		\]
		so $\omega(T)=1$, but for every unit vector $x$, $|\langle Tx,x\rangle|<1$. Thus $T\notin\theta(\ell^2)$, showing that infinite orthogonal direct sums need not preserve the property.
	\end{example}
	
	\begin{remark}\label{rem:finite-direct-sum-normal}
		For normal $T_1,T_2$, the conclusion of Theorem~\ref{thm:finite-direct-sum}
		also follows directly from the spectral characterization
		Theorem~\ref{thm:normal}: both operators have pure point spectrum, their
		eigenvalue-modulus sets are reverse well-ordered, and the finite union of two
		reverse well-ordered subsets of $\mathbb{R}$ is again reverse well-ordered.
		The point of Theorem~\ref{thm:finite-direct-sum} is that normality is not
		needed.
	\end{remark}

	\section{Essential numerical radius and compact perturbations}\label{sec:compact}

	In this section we study compact perturbations and their effect on membership
	in $\theta(H)$. Throughout the section we assume that
	\[
	T\in\theta(H),\qquad K\in\mathcal K(H),\qquad S=T+K.
	\]
	Thus $T$ is a fixed member of $\theta(H)$, and we ask when the compact
	perturbation $S$ remains in $\theta(H)$. The central notion is
	the essential numerical radius, which is the numerical-radius analogue of
	the essential numerical range. After obtaining the criterion for
	$S\in\theta(H)$, we prove that the compact-perturbation-invariant core
	$\mathfrak C_\theta(H)$ is exactly $\mathcal K(H)$, without assuming
	normality.
	
	If $H$ is finite-dimensional, then $\mathcal K(H)=\mathcal B(H)$ and
	$\theta(H)=\mathcal B(H)$. Thus the compact-perturbation criterion has no
	infinite-dimensional reducing subspaces to test, and
	$\mathfrak C_\theta(H)=\mathcal B(H)=\mathcal K(H)$. The essential
	numerical range and the essential numerical radius are used below only when
	$H$ is infinite-dimensional, and all arguments involving weakly null unit
	sequences are carried out in that setting. Unless the finite-dimensional case
	is explicitly mentioned, we assume this from now on in the present section.
	
	For $A\in\mathcal{B}(H)$, the \emph{essential norm} of $A$ is defined by
	\[
	\|A\|_e
	:=
	\inf\{\|A+K\|:K\in\mathcal{K}(H)\}.
	\]
	Equivalently, if
	\[
	\pi:\mathcal{B}(H)\longrightarrow
	\mathcal{B}(H)/\mathcal{K}(H)
	\]
	denotes the quotient map onto the Calkin algebra, then
	\[
	\|A\|_e=\|\pi(A)\|.
	\]
	In particular, $\|A\|_e=0$ if and only if $A$ is compact.
	
	Fillmore, Stampfli and Williams \cite{FillmoreStampfliWilliams1972} introduced the notion of the \emph{essential numerical range}. For $A\in\mathcal{B}(H)$, it is defined by
	\[
	W_e(A)
	:=
	\bigcap_{K\in\mathcal{K}(H)}
	\overline{W(A+K)}.
	\]
	Equivalently,
	\[
	W_e(A)
	=
	\left\{
	z\in\mathbb{C}:
	\begin{array}{l}
		\text{there exists a sequence of unit vectors }x_n\rightharpoonup0\\
		\text{such that }\langle Ax_n,x_n\rangle\to z
	\end{array}
	\right\}.
	\]
	As shown in \cite{FillmoreStampfliWilliams1972} (see also \cite{Mansoori2017}), $W_e(A)$ is a nonempty compact convex subset of $\mathbb{C}$.
	
	The basic properties of $W_e(A)$---compactness, convexity, and invariance
	under compact perturbations---are standard; see
	\cite{FillmoreStampfliWilliams1972,Bourin2003}. The compact-perturbation
	characterization recalled in \cite[Proposition~1.1]{Bourin2003} states
	that for every $A$ there is a compact $K$ with
	$W_e(A)=\overline{W(A+K)}$.
	
	For $A\in\mathcal{B}(H)$, we define the \emph{essential numerical radius} by
	\[
	\omega_e(A)
	:=
	\max\{|z|:z\in W_e(A)\}.
	\]
	The maximum exists because $W_e(A)$ is compact. Thus
	\[
	0\leqslant\omega_e(A)\leqslant\omega(A).
	\]
	
	The following variational formula relates the definition through the essential numerical range to compact perturbations. This formula is standard; we include a proof for completeness.
	
	\begin{lemma}\label{lem:essential-radius-inf}
		For every $A\in\mathcal B(H)$,
		\[
		\omega_e(A)=\inf_{K\in\mathcal K(H)}\omega(A+K).
		\]
	\end{lemma}
	
	\begin{proof}
		For every compact $K$, the definition
		$W_e(A)=\bigcap_{L\in\mathcal K(H)}\overline{W(A+L)}$ gives
		$W_e(A)\subseteq\overline{W(A+K)}$. Hence
		$\omega_e(A)\leqslant\omega(A+K)$, and therefore
		\[
		\omega_e(A)\leqslant\inf_{K\in\mathcal K(H)}\omega(A+K).
		\]
		
		For the reverse inequality, put $r=\omega_e(A)$ and fix $\varepsilon>0$.
		There exists a finite-dimensional subspace $F\subset H$ such that
		\[
		\omega\bigl((P_{F^\perp}A)|_{F^\perp}\bigr)\leqslant r+\varepsilon,
		\]
		where the displayed operator is regarded as an operator on $F^\perp$.
		Indeed, if no such $F$ existed, we could recursively choose an orthonormal
		sequence $(x_n)$ satisfying
		$|\langle Ax_n,x_n\rangle|>r+\varepsilon$ for every $n$:
		when $x_1,\ldots,x_{n-1}$ have been chosen, apply the assumed failure to
		$F_{n-1}=\operatorname{span}\{x_1,\ldots,x_{n-1}\}$. The resulting
		orthonormal sequence is weakly null. Since the scalar sequence
		$\langle Ax_n,x_n\rangle$ is bounded, a subsequence converges to some $z$.
		Then $z\in W_e(A)$ by the weakly-null characterization of the essential
		numerical range, while $|z|\geqslant r+\varepsilon$, contradicting the
		definition of $r$.
		
		Let $P=P_F$ and $Q=I-P$. The operator $K_F=-PA-QAP$ has finite rank, and
		\[
		A+K_F=QAQ=0_F\oplus (QA)|_{F^\perp}.
		\]
		Consequently,
		\[
		\omega(A+K_F)=\omega\bigl((QA)|_{F^\perp}\bigr)\leqslant r+\varepsilon.
		\]
		As $\varepsilon>0$ was arbitrary, the reverse inequality follows.
	\end{proof}
	
	For $\zeta\in\mathbb{T}$, recall that
	\[
	H_\zeta(A)
	=
	\frac{\zeta A+\overline{\zeta}A^*}{2}.
	\]
	Since $H_\zeta(A)$ is self-adjoint, its Fredholm essential spectrum
	$\sigma_e(H_\zeta(A))$ is a compact subset of $\mathbb{R}$. The essential
	numerical range satisfies the standard support-function identity
	\[
	\sup\sigma_e(H_\zeta(A))
	=
	\max_{z\in W_e(A)}
	\operatorname{Re}(\zeta z),
	\qquad \zeta\in\mathbb{T},
	\]
	see~\cite{FillmoreStampfliWilliams1972,Mansoori2017}. Consequently,
	\[
	\begin{aligned}
		\max_{\zeta\in\mathbb{T}}
		\sup\sigma_e(H_\zeta(A))
		&=
		\max_{\zeta\in\mathbb{T}}
		\max_{z\in W_e(A)}
		\operatorname{Re}(\zeta z)\\&=\max_{z\in W_e(A)}
		\max_{\zeta\in\mathbb{T}}
		\operatorname{Re}(\zeta z)=
		\max_{z\in W_e(A)}|z|=\omega_e(A).
	\end{aligned}
	\]
	Hence we have the equivalent formula
	\[
	\omega_e(A)
	=\max_{\zeta\in\mathbb{T}}
	\sup\sigma_e(H_\zeta(A))=
	\max_{\zeta\in\mathbb{T}}
	\sup\sigma_e\left(
	\frac{\zeta A+\overline{\zeta}A^*}{2}
	\right).
	\]
	
	The preceding formula also shows that $\omega_e$ is invariant under compact
	perturbations. More precisely, if $A,B\in\mathcal{B}(H)$ and $A-B$ is compact,
	then
	\[
	\omega_e(A)=\omega_e(B).
	\]
	Indeed, for every $\zeta\in\mathbb{T}$,
	\[
	H_\zeta(A)-H_\zeta(B)
	=
	\frac{\zeta(A-B)+\overline{\zeta}(A-B)^*}{2}
	\]
	is compact. Therefore
	\[
	\sigma_e(H_\zeta(A))
	=
	\sigma_e(H_\zeta(B)),
	\]
	by invariance of the Fredholm essential spectrum under compact perturbations.
	Taking the maximum over $\zeta\in\mathbb{T}$ gives
	\[
	\omega_e(A)=\omega_e(B).
	\]
	
	The following proposition is the key observation relating the numerical radius and the essential numerical radius to attainment. It says that a strict gap between the two radii always forces the numerical radius to be attained.
	
	\begin{proposition}\label{prop:strict-gap-attainment}
		If $A\in\mathcal{B}(H)$ satisfies
		\[
		\omega(A)>\omega_e(A),
		\]
		then $A$ attains its numerical radius.
	\end{proposition}
	
	\begin{proof}
		By Lemma~\ref{lem:numerical}, there exists $\zeta_0\in\mathbb{T}$ such that
		\[
		\omega(A)=\sup\sigma\bigl(H_{\zeta_0}(A)\bigr).
		\]
		The support-function formula for the essential numerical radius gives
		\[
		\sup\sigma_e\bigl(H_{\zeta_0}(A)\bigr)
		\leqslant\omega_e(A)<\omega(A).
		\]
		Thus the top spectral point $\omega(A)$ of the self-adjoint operator
		$H_{\zeta_0}(A)$ lies outside its essential spectrum. It is therefore an
		isolated eigenvalue of finite multiplicity. Lemma~\ref{lem:numerical} now
		shows that $\omega(A)$ is attained.
	\end{proof}
	
	\emph{Convention.}
	If $M$ is a finite-dimensional nonzero subspace, then for every
	$A\in\mathcal{B}(M)$ the numerical radius $\omega(A)$ is attained on the
	unit sphere of $M$. Hence finite-dimensional reducing subspaces impose no
	additional condition in the definition of $\theta(H)$. In the statements
	below, whenever $M$ is infinite-dimensional, $P_M$ denotes the orthogonal
	projection onto $M$, and
	\[
	P_M T|_M:M\to M
	\]
	denotes the compression of $T$ to $M$. Notice that $M$ need not reduce $T$.
	
	By Theorem~\ref{thm:general}, membership in $\theta(H)$ is characterized by
	the attainment of $\omega(S|_M)$ on every nonzero
	$M\in\operatorname{Red}(S)$. The next proposition shows that, under a
	compact perturbation, the only possible obstruction occurs at the level
	of the essential numerical radius.
	
	\begin{proposition}\label{prop:compactperturbation}
		Let $T\in\theta(H)$, $K\in\mathcal{K}(H)$, and put
		\[
		S=T+K.
		\]
		For every nonzero $M\in\operatorname{Red}(S)$, if $\dim M<\infty$, then
		$\omega(S|_M)$ is attained. If $\dim M=\infty$, define
		\[
		C_M=P_MT|_M,
		\qquad
		K_M=P_MK|_M.
		\]
		Then $K_M$ is compact,
		\[
		S|_M=C_M+K_M,
		\]
		and
		\[
		\omega_e(S|_M)=\omega_e(C_M).
		\]
		Moreover, if
		\[
		\omega(S|_M)>\omega_e(C_M),
		\]
		then $\omega(S|_M)$ is attained on $M$.
	\end{proposition}
	
	\begin{proof}
		If $\dim M<\infty$, then the numerical radius is attained on the compact
		unit sphere of $M$. Now assume that $\dim M=\infty$. Since
		$M\in\operatorname{Red}(S)$, we have
		\[
		S(M)\subset M.
		\]
		Therefore
		\[
		S|_M
		=
		P_MS|_M
		=
		P_M(T+K)|_M
		=
		C_M+K_M.
		\]
		Because $K$ is compact, $K_M$ is compact. By the compact-perturbation
		invariance of the essential numerical radius proved above,
		\[
		\omega_e(S|_M)
		=
		\omega_e(C_M).
		\]
		
		Suppose that
		\[
		r:=\omega(S|_M)>\omega_e(C_M).
		\]
		By Lemma~\ref{lem:numerical}, there exists
		$\zeta_0\in\mathbb{T}$ such that
		\[
		r
		=
		\sup\sigma\left(H_{\zeta_0}(S|_M)\right).
		\]
		Moreover,
		\[
		\sup\sigma_e\left(H_{\zeta_0}(S|_M)\right)
		\leqslant
		\omega_e(S|_M)
		=
		\omega_e(C_M)
		<
		r.
		\]
		Since $H_{\zeta_0}(S|_M)$ is self-adjoint, every spectral point strictly
		above the essential spectrum is an isolated eigenvalue of finite
		multiplicity. Hence
		\[
		r\in
		\sigma_p\left(H_{\zeta_0}(S|_M)\right).
		\]
		Lemma~\ref{lem:numerical} then implies that $\omega(S|_M)$ is attained.
	\end{proof}
	
	The preceding proposition immediately yields the following criterion for membership in $\theta(H)$ under a compact perturbation.
	
	\begin{theorem}\label{thm:compactcriterion}
		Let $T\in\theta(H)$, $K\in\mathcal{K}(H)$, and put
		\[
		S=T+K.
		\]
		Then $S\in\theta(H)$ if and only if, for every infinite-dimensional
		nonzero $M\in\operatorname{Red}(S)$ satisfying
		\[
		\omega(S|_M)=\omega_e(P_MT|_M),
		\]
		there exists $\zeta\in\mathbb{T}$ such that
		\[
		\omega(S|_M)
		\in
		\sigma_p\left(H_\zeta(S|_M)\right).
		\]
		Finite-dimensional nonzero reducing subspaces impose no additional condition.
	\end{theorem}
	
	\begin{proof}
		If $M$ is finite-dimensional, then $\omega(S|_M)$ is automatically attained.
		
		Now let $M\in\operatorname{Red}(S)$ be infinite-dimensional. By
		Proposition~\ref{prop:compactperturbation},
		\[
		\omega(S|_M)
		\geqslant
		\omega_e(S|_M)
		=
		\omega_e(P_MT|_M).
		\]
		If the inequality is strict, then
		Proposition~\ref{prop:compactperturbation} shows that
		$\omega(S|_M)$ is automatically attained. Thus the only possible obstruction
		to membership in $\theta(H)$ occurs when
		\[
		\omega(S|_M)=\omega_e(P_MT|_M).
		\]
		For such a reducing subspace $M$, Lemma~\ref{lem:numerical} says that
		$\omega(S|_M)$ is attained if and only if there exists
		$\zeta\in\mathbb{T}$ such that
		\[
		\omega(S|_M)
		\in
		\sigma_p\left(H_\zeta(S|_M)\right).
		\]
		Applying this condition to every infinite-dimensional nonzero
		$M\in\operatorname{Red}(S)$ proves the theorem.
	\end{proof}
	
	\begin{example}\label{ex:compact-criterion}
		Let $H=H_1\oplus H_2$ with $H_1=H_2=\ell^2(\mathbb{N})$, let $\{e_n\}$ and $\{f_n\}$ be the standard orthonormal bases of $H_1$ and $H_2$, respectively, and let
		\[
		T=I_{H_1}\oplus \operatorname{diag}\left(1,\frac12,\frac13,\dots\right).
		\]
		Then $T$ is normal, has pure point spectrum, and its eigenvalue-modulus set is $\{1,\frac12,\frac13,\dots\}$, which is reverse well-ordered. By Theorem~\ref{thm:normal}, $T\in\theta(H)$. Let $K=-P_{e_1}\oplus0$. Then $K$ has rank one and is compact. Put $S=T+K$. Thus
		\[
		S=(I_{H_1}-P_{e_1})\oplus \operatorname{diag}\left(1,\frac12,\frac13,\dots\right).
		\]
		We verify $S\in\theta(H)$ using Theorem~\ref{thm:compactcriterion}. Let $M\in\operatorname{Red}(S)$ be infinite-dimensional. By Lemma~\ref{lem:reducing}, $M$ is the orthogonal direct sum of its intersections with the eigenspaces of $S$. Every eigenvalue of $S$ other than $1$ has a one-dimensional eigenspace. We distinguish two cases.
		
		\begin{enumerate}
			\item $M$ contains a unit vector $x$ with $Sx=x$. Then $\omega(S|_M)=1$, and this value is attained at $x$. In particular, if $\omega(S|_M)=\omega_e(P_MT|_M)$, then the eigenvalue condition of Theorem~\ref{thm:compactcriterion} holds: since $S$ is normal and $Sx=x$, we have $H_{\zeta=1}(S|_M)x=x$; if the inequality is strict, Proposition~\ref{prop:compactperturbation} already gives attainment.
			\item $M$ contains no eigenvector of $S$ with eigenvalue $1$. Then $M$ is contained in the closed span of the eigenvectors with eigenvalues $0,\frac12,\frac13,\dots$. Hence the $H_1$-component of $M$ is contained in $\mathbb{C}e_1$, and the $H_2$-component is contained in $\overline{\operatorname{span}}\{f_2,f_3,\dots\}$, on which $T$ is the compact diagonal operator $\operatorname{diag}(\frac12,\frac13,\dots)$. Thus $P_MT|_M$ is compact, so $\omega_e(P_MT|_M)=0$. Since $M$ is infinite-dimensional, $S|_M\neq0$, hence $\omega(S|_M)>0$, and the strict-gap case again gives attainment.
		\end{enumerate}
		
		Thus the condition of Theorem~\ref{thm:compactcriterion} holds for every infinite-dimensional $M\in\operatorname{Red}(S)$. Finite-dimensional reducing subspaces impose no condition, so $S\in\theta(H)$.
	\end{example}
	
	We now turn to the compact-perturbation-invariant core $\mathfrak{C}_\theta(H)$ of
	$\theta(H)$. The next theorem proves that this core consists exactly of the
	compact operators, with no normality assumption. The proof uses a finite-rank
	projection lemma with prescribed tail bounds and a common-majorant lemma for
	norm-compact families of positive compact operators.
	Since $K=0$ is allowed, the condition $T+K\in\theta(H)$ for every compact
	$K$ already forces $T\in\theta(H)$; we nevertheless include this restriction
	explicitly in the definition, to emphasize that the core is a subset of
	$\theta(H)$.
	
	\begin{definition}
		The \emph{compact-perturbation-invariant core} of $\theta(H)$ is
		\[
		\mathfrak{C}_\theta(H)
		:=
		\left\{T\in\theta(H):
		T+K\in\theta(H)\text{ for every }K\in\mathcal{K}(H)
		\right\}.
		\]
	\end{definition}
	
	The following lemma supplies the projection sequence used in the proof of
	the common-majorant lemma.
	
	\begin{lemma}[Finite-rank projections with prescribed tail bounds]\label{lem:tail-projections}
		Let $H$ be an infinite-dimensional separable complex Hilbert space and let
		$\mathcal F$ be a nonempty norm-compact family of positive compact operators
		on $H$. Then there exists
		a sequence of finite-rank orthogonal projections $\{P_n\}_{n=0}^\infty$
		such that
		\begin{enumerate}
			\item $P_0=0$ and $P_{n-1}\leqslant P_n$ for every $n\geqslant1$;
			\item $P_n\xrightarrow{\mathrm s}I$;
			\item for every $n\geqslant1$,
			\[
			\sup_{A\in\mathcal F}\|(I-P_n)A^{\frac12}\|\leqslant2^{-n^2}.
			\]
		\end{enumerate}
	\end{lemma}
	
	\begin{proof}
		Fix an orthonormal basis $\{e_k\}_{k=1}^\infty$ of $H$ and put
		\[
		\mathcal G=\{A^{\frac12}:A\in\mathcal F\}.
		\]
		The square-root map is norm continuous on the positive operators, so
		$\mathcal G$ is norm compact; moreover, every member of $\mathcal G$ is
		compact. We construct the projections inductively. Set $P_0=0$. Suppose
		that $P_{n-1}$ has been constructed for some $n\geqslant1$. Since
		$\mathcal G$ is norm compact, choose finitely many operators
		$G_1,\ldots,G_m\in\mathcal G$ such that for every $G\in\mathcal G$ there
		is some $j\in\{1,\ldots,m\}$ with
		\[
		\|G-G_j\|\leqslant2^{-n^2-1}.
		\]
		For each $j$, the compactness of $G_j$ gives a finite-rank orthogonal
		projection $Q_{n,j}$ such that
		\[
		\|(I-Q_{n,j})G_j\|\leqslant2^{-n^2-1}.
		\]
		Let $P_n$ be the orthogonal projection onto
		\[
		\left(\sum_{j=1}^m\operatorname{ran}Q_{n,j}\right)
		+P_{n-1}H+\operatorname{span}\{e_1,\ldots,e_n\}.
		\]
		Then $P_n$ has finite rank and $P_{n-1}\leqslant P_n$. Since
		$Q_{n,j}\leqslant P_n$, we have $I-P_n\leqslant I-Q_{n,j}$ and hence
		$(I-P_n)=(I-P_n)(I-Q_{n,j})$. Therefore
		\[
		\|(I-P_n)G_j\|
		\leqslant\|(I-Q_{n,j})G_j\|
		\leqslant2^{-n^2-1}
		\qquad(j=1,\ldots,m).
		\]
		
		This completes the inductive construction. Property (1) holds by
		construction. For (2), let $x\in H$ and $\varepsilon>0$. Choose $N$ such
		that
		\[
		\left\|x-\sum_{k=1}^N\langle x,e_k\rangle e_k\right\|<\frac{\varepsilon}{2}.
		\]
		Since $e_1,\ldots,e_N\in\operatorname{ran}P_n$ for every $n\geqslant N$,
		the projection $P_n$ fixes the partial sum; hence
		$\|P_nx-x\|\leqslant2\|x-\sum_{k=1}^N\langle x,e_k\rangle e_k\|<\varepsilon$
		for all $n\geqslant N$. Thus
		$P_n\xrightarrow{\mathrm s}I$. For (3), fix $n\geqslant1$ and
		$A\in\mathcal F$, and put $G=A^{\frac12}\in\mathcal G$. Choose
		$j\in\{1,\ldots,m\}$ with $\|G-G_j\|\leqslant2^{-n^2-1}$. Then
		\[
		\begin{aligned}
			\|(I-P_n)G\|
			&\leqslant\|(I-P_n)(G-G_j)\|+\|(I-P_n)G_j\|\\
			&\leqslant\|G-G_j\|+\|(I-Q_{n,j})G_j\|\\
			&\leqslant2^{-n^2-1}+2^{-n^2-1}
			=2^{-n^2}.
		\end{aligned}
		\]
		Taking the supremum over $A\in\mathcal F$ proves (3).
	\end{proof}
	
	\begin{lemma}[A common positive compact majorant]\label{lem:compact-majorant}
		Let $\mathcal F$ be a nonempty norm-compact family of positive compact operators on an
		infinite-dimensional separable Hilbert space $H$. Then there is an injective
		positive compact
		operator $B$ such that
		\[
		0\leqslant A\leqslant B\qquad(A\in\mathcal F).
		\]
	\end{lemma}
	
	\begin{proof}
		By Lemma~\ref{lem:tail-projections}, applied to $\mathcal F$, there is an
		increasing sequence of finite-rank orthogonal projections
		\[
		P_0=0,\qquad P_1\leqslant P_2\leqslant\cdots,
		\qquad P_n\xrightarrow{\mathrm s}I,
		\]
		satisfying
		\[
		\sup_{A\in\mathcal F}\|(I-P_n)A^{\frac12}\|\leqslant2^{-n^2}
		\qquad(n\geqslant1).
		\]
		
		Set $Q_n=P_n-P_{n-1}$. The $Q_n$ are mutually orthogonal finite-rank
		projections and $\sum_nQ_n=I$ in the strong operator topology. Put
		\[
		c_n=\sup_{A\in\mathcal F}\|A^{\frac12}Q_n\|.
		\]
		The number $c_1$ is finite because $\mathcal F$ is norm compact. For
		$n\geqslant2$, self-adjointness of $A^{\frac12}$ and $Q_n\leqslant I-P_{n-1}$ give
		\[
		\begin{aligned}
			c_n&=\sup_{A\in\mathcal F}\|Q_nA^{\frac12}\|\leqslant\sup_{A\in\mathcal F}\|(I-P_{n-1})A^{\frac12}\|\leqslant2^{-(n-1)^2}.
		\end{aligned}
		\]
		Thus, with $d_n=2^{-n}$,
		\[
		\sum_{n=1}^{\infty}\frac{c_n^2}{d_n}<\infty.
		\]
		
		Define
		\[
		D=\sum_{n=1}^{\infty}d_nQ_n.
		\]
		Since the $Q_n$ are orthogonal finite-rank projections and $d_n\to0$,
		$D$ is positive and compact. Also,
		\[
		\langle Dx,x\rangle=\sum_{n=1}^{\infty}d_n\|Q_nx\|^2>0
		\qquad(x\neq0),
		\]
		because $\sum_nQ_n=I$ strongly. Hence $D$ is injective.
		
		For $A\in\mathcal F$ and $x\in H$, the strong expansion
		$x=\sum_nQ_nx$ and the Cauchy--Schwarz inequality yield
		\[
		\begin{aligned}
			\|A^{\frac12}x\|
			&\leqslant\sum_{n=1}^{\infty}c_n\|Q_nx\|\\
			&\leqslant
			\left(\sum_{n=1}^{\infty}\frac{c_n^2}{d_n}\right)^{\frac12}
			\left(\sum_{n=1}^{\infty}d_n\|Q_nx\|^2\right)^{\frac12}\\
			&=L^{\frac12}\langle Dx,x\rangle^{\frac12},
		\end{aligned}
		\]
		where $L=\sum_n \frac{c_n^2}{d_n}<\infty$. Squaring gives
		$\langle Ax,x\rangle\leqslant L\langle Dx,x\rangle$ for every $x$, so
		$A\leqslant LD$. Taking $B=(L+1)D$ gives an injective positive compact
		operator that dominates every $A\in\mathcal F$, including the case
		$L=0$.
	\end{proof}
	
	\begin{theorem}[The compact-perturbation-invariant core]\label{thm:invariant}
		Let $H$ be a separable complex Hilbert space. Then
		\[
		\mathfrak{C}_\theta(H)=\mathcal K(H).
		\]
	\end{theorem}
	
	\begin{proof}
		If $H$ is finite-dimensional, then $\mathfrak C_\theta(H)=\mathcal B(H)=\mathcal K(H)$, so the result is immediate. Assume henceforth that $H$ is infinite-dimensional. If $T$ is compact, then $T+K$ is compact for every compact $K$.
		Proposition~\ref{prop:compact} therefore gives $T+K\in\theta(H)$ for every
		$K\in\mathcal K(H)$. Hence
		$\mathcal K(H)\subseteq\mathfrak C_\theta(H)$.
		
		For the converse, let $T\in\theta(H)\setminus\mathcal K(H)$ and put
		\[
		r=\omega_e(T)=\inf_{K\in\mathcal K(H)}\omega(T+K),
		\]
		where the identity is Lemma~\ref{lem:essential-radius-inf}. If
		$\|T\|_e:=\inf_{K\in\mathcal K(H)}\|T+K\|$ denotes the essential norm,
		then
		\[
		r\geqslant\frac12\|T\|_e>0,
		\]
		because $T$ is not compact.
		
		For $\zeta\in\mathbb{T}$, put
		\[
		H_\zeta=\operatorname{Re}(\zeta T)
		=\frac{\zeta T+\overline\zeta T^*}{2}.
		\]
		For every compact $K$,
		\[
		\begin{aligned}
			\|H_\zeta\|_e
			\leqslant\|\operatorname{Re}(\zeta(T+K))\|\leqslant\omega(T+K).
		\end{aligned}
		\]
		Taking the infimum over $K$ gives $\|H_\zeta\|_e\leqslant r$. Let
		$\pi:\mathcal B(H)\to\mathcal B(H)/\mathcal K(H)$ be the quotient map.
		Since $H_\zeta$ is self-adjoint,
		$\sigma(\pi(H_\zeta))\subseteq[-r,r]$. Therefore continuous functional
		calculus gives
		\[
		P_\zeta:=(H_\zeta-rI)_+\in\mathcal K(H).
		\]
		Here the subscript $+$ denotes the positive part of a self-adjoint operator.
		The map $\zeta\mapsto P_\zeta$ is norm continuous and its domain,
		the unit circle $\mathbb{T}$, is compact. Its image
		$\mathcal F=\{P_\zeta:\zeta\in\mathbb{T}\}$ is therefore a norm-compact
		family of positive compact operators. Lemma~\ref{lem:compact-majorant}
		gives an injective positive compact operator $B$ such that
		\[
		P_\zeta\leqslant B\qquad(\zeta\in\mathbb{T}).
		\]
		As $H_\zeta-rI\leqslant(H_\zeta-rI)_+=P_\zeta$, we obtain
		\[
		rI+2B-H_\zeta\geqslant B\qquad(\zeta\in\mathbb{T}).
		\]
		
		Set $C=2B$ and define
		\[
		X=\left(I+\frac Cr\right)^{-\frac12}.
		\]
		The operator $X$ is positive and invertible. Moreover, $X-I$ is compact,
		since the continuous function $t\mapsto\left(1+\frac tr\right)^{-\frac12}-1$ vanishes at $0$
		and $C$ is compact. We also have
		\[
		X(rI+C)X=rI.
		\]
		Define $S=XTX$. Then
		\[
		S-T=(X-I)TX+T(X-I)\in\mathcal K(H),
		\]
		so $S$ is a compact perturbation of $T$. For every $\zeta\in\mathbb{T}$,
		\[
		\begin{aligned}
			rI-\operatorname{Re}(\zeta S)
			&=X(rI+C-H_\zeta)X\\
			&\geqslant XBX.
		\end{aligned}
		\]
		Because $B$ is positive and injective and $X$ is invertible,
		$\langle XBXx,x\rangle>0$ for every $x\neq0$. Thus for every unit vector
		$x$ and every $\zeta\in\mathbb{T}$,
		\[
		\operatorname{Re}\bigl(\zeta\langle Sx,x\rangle\bigr)
		\leqslant r-\langle XBXx,x\rangle<r.
		\]
		Taking the maximum over $\zeta\in\mathbb{T}$ gives
		\[
		|\langle Sx,x\rangle|<r\qquad(\|x\|=1).
		\tag{*}
		\]
		Consequently $\omega(S)\leqslant r$. On the other hand, $S-T$ is compact,
		so $\omega_e(S)=\omega_e(T)=r$, and therefore $\omega(S)\geqslant r$.
		It follows that $\omega(S)=r$, but (*) shows that the supremum is not
		attained on the unit sphere. Since $H$ itself is a nonzero reducing subspace
		of $S$, we have $S\notin\theta(H)$. As $S-T$ is compact, this proves
		$T\notin\mathfrak C_\theta(H)$.
	\end{proof}
	
	In particular, the normal operators in $\mathfrak C_\theta(H)$ are
	exactly the compact normal operators.
	
	\begin{remark}[Related work]\label{rem:related-work}
		The compact-perturbation questions considered here are related in spirit to
		the study of operators with closed numerical ranges \cite{JiLiang2018} and
		to the study of numerical ranges in semi-Hilbertian spaces
		\cite{SenBirbonshiPaul2024}. We do not identify these properties with the
		one studied here. Theorem~\ref{thm:invariant} determines the operators
		all of whose compact perturbations belong to $\theta(H)$, and shows that this
		compact-perturbation-invariant core is precisely $\mathcal K(H)$.
	\end{remark}
	
	The following proposition gives a non-normal, noncompact operator in
	$\theta(H)$ whose membership is destroyed by a rank-one perturbation.
	Scaling both the operator and its perturbation gives such examples with
	arbitrarily small perturbation norms.

	\begin{proposition}[Failure of rank-one stability for a non-normal operator]
		\label{prop:rankone-instability}
		Let $H_0=\ell^2(\mathbb{Z})$ with its standard orthonormal basis
		$\{e_n\}_{n\in\mathbb{Z}}$, and let $Ue_n=e_{n+1}$ be the bilateral shift.
		For $x=\sum_{n\in\mathbb{Z}}x_ne_n$, write $\varphi(x)=x_0$, and set
		$H=H_0\oplus\mathbb{C}$. Define
		\[
		T(x,\alpha)=(Ux,\varphi(x)+2\alpha),
		\qquad
		F(x,\alpha)=(0,-\varphi(x)).
		\]
		Then $F$ has rank one,
		\[
		T+F=U\oplus 2I_{\mathbb{C}},
		\]
		and $T\in\theta(H)$ but $T+F\notin\theta(H)$.
	\end{proposition}
	
	\begin{proof}
		Let $f=(0,1)$. If $T(x,\alpha)=\lambda(x,\alpha)$, then
		$Ux=\lambda x$. Since the bilateral shift has no eigenvalues, $x=0$;
		for a nonzero eigenvector this forces $\lambda=2$. Thus
		$\sigma_p(T)=\{2\}$ and $\ker(T-2I)=\mathbb{C}f$, so $T$ does not
		have pure point spectrum. Also, $T^*f=(e_0,2)$ and
		$T(e_0,0)=(e_1,1)$. Moreover, $\|T(e_0,0)\|=\sqrt{2}$, whereas
		$\|T^*(e_0,0)\|=1$. Recall that a normal operator satisfies
		$\|Tx\|=\|T^*x\|$ for every $x$, since
		$T^*T=TT^*$; hence $T$ is non-normal. The operator $T$ is also non-compact: its compression to $H_0\oplus\{0\}$ is the non-compact bilateral shift. Indeed, a compression of a compact operator to a closed subspace is compact, so this compression rules out compactness of $T$.
		
		We first show that $T$ is irreducible, meaning that its only reducing
		subspaces are $\{0\}$ and $H$. Let $M$ be a reducing subspace and let
		$P$ be the orthogonal projection onto $M$. Then $P$ commutes with $T$,
		so $Pf\in\ker(T-2I)=\mathbb{C}f$. Since $P^2=P$, either $Pf=0$ or
		$Pf=f$.
		
		Suppose first that $Pf=f$. Then $f\in M$. Since $M$ reduces $T$, it
		also contains $T^*f-2f=(e_0,0)$. It follows from
		$T(e_0,0)=(e_1,0)+f$ that $(e_1,0)\in M$, and from
		$T^*(e_0,0)=(e_{-1},0)$ that $(e_{-1},0)\in M$. Iterating $T$ and
		$T^*$ shows that $(e_n,0)\in M$ for every $n\in\mathbb{Z}$; hence
		$M=H$.
		
		Suppose instead that $Pf=0$. Then $M\subseteq H_0\oplus\{0\}$. If
		$(x,0)\in M$, invariance under $T$ gives
		$(Ux,\varphi(x))\in M$, so $\varphi(x)=0$ and $(Ux,0)\in M$. Invariance
		under $T^*$ also gives $(U^*x,0)\in M$. Thus $(U^nx,0)\in M$ and
		$\varphi(U^nx)=0$ for every
		$n\in\mathbb{Z}$. Since $\varphi(U^nx)$ runs through all coordinates of
		$x$, we get $x=0$. Hence $M=\{0\}$. This proves irreducibility.
		
		Next, put $S=U\oplus 2I_{\mathbb{C}}$. Since $T-S=-F$ is finite rank, compact-perturbation invariance of the essential numerical radius gives
		$\omega_e(T)=\omega_e(S)=1$. Indeed, the unit vectors $x_N=N^{-\frac12}\sum_{k=N+1}^{2N} e_k$ are weakly null and satisfy $\langle Ux_N,x_N\rangle=\frac{N-1}N$, so $\omega_e(U)=1$. A weakly null sequence in $H_0\oplus\mathbb{C}$ has scalar coordinate tending to zero; after normalizing its $H_0$-coordinates, the definition of $W_e$ shows that $W_e(U\oplus 2I_{\mathbb C})=W_e(U)$. Thus $\omega_e(S)=\omega_e(U)=1$. Since $Tf=2f$, we have $\omega(T)\geqslant2>\omega_e(T)$.
		Proposition~\ref{prop:strict-gap-attainment} implies that $T$ attains its
		numerical radius. Because $T$ is irreducible, $H$ is its only nonzero
		reducing subspace. Consequently, $T\in\theta(H)$.
		
		On the other hand, $M=H_0\oplus\{0\}$ reduces $S=T+F$, and
		$S|_M=U$. We have $\omega(U)=1$, but this value is not attained. In
		fact, if a unit vector $x$ satisfied $|\langle Ux,x\rangle|=1$, equality
		in the Cauchy--Schwarz inequality would imply $Ux=\lambda x$ for some
		$|\lambda|=1$, contradicting that $U$ has no eigenvalues. Hence
		$S\notin\theta(H)$. Moreover, for every $t>0$, the scalar-multiplication property gives $tT\in\theta(H)$, while $tT+tF=t(T+F)$ is not in $\theta(H)$ because $H_0\oplus\{0\}$ still reduces $t(T+F)$ and its restriction is $tU$, whose numerical radius is not attained. Since $\|tF\|=t\|F\|\to0$ as $t\downarrow0$, this gives rank-one perturbations of arbitrarily small norm that destroy membership in $\theta(H)$.
	\end{proof}
	
	\begin{example}[Failure of power stability without normality]
		\label{ex:nonnormal-power}
		Let $T$ be the non-normal irreducible operator constructed in
		Proposition~\ref{prop:rankone-instability}; thus $T\in\theta(H)$. We show
		that
		\[
		T^2\notin\theta(H).
		\]
		Write $H_0=\ell^2(\mathbb Z)$ and split it into its even and odd coordinate
		subspaces,
		\[
		H_{\mathrm e}=\overline{\operatorname{span}}\{e_{2k}:k\in\mathbb Z\},
		\qquad
		H_{\mathrm o}=\overline{\operatorname{span}}\{e_{2k+1}:k\in\mathbb Z\}.
		\]
		Let $V_{\mathrm e}=U^2|_{H_{\mathrm e}}$ and
		$V_{\mathrm o}=U^2|_{H_{\mathrm o}}$. Under the natural identifications
		$e_{2k}\leftrightarrow e_k$ and $e_{2k+1}\leftrightarrow e_k$, both are
		bilateral shifts. Define a bounded operator $R:H_{\mathrm e}\to H_{\mathrm o}$
		by
		\[
		Re_{2k}=-2e_{2k-1},\qquad k\in\mathbb Z.
		\]
		Then $\|R\|=2$ and
		\[
		V_{\mathrm o}R=RV_{\mathrm e},
		\qquad
		V_{\mathrm o}^*R=RV_{\mathrm e}^*.
		\]
		Consequently the graph
		\[
		M_0=\{x+Rx:x\in H_{\mathrm e}\}\subset H_0
		\]
		is a closed reducing subspace for $U^2$.
		
		A direct computation gives
		\[
		T^2(x,\alpha)
		=\bigl(U^2x,\,\varphi(Ux)+2\varphi(x)+4\alpha\bigr).
		\]
		Let $\psi(x)=\varphi(Ux)+2\varphi(x)$. If
		$x=\sum_k a_ke_{2k}\in H_{\mathrm e}$, then
		\[
		\varphi(x+Rx)=a_0,
		\qquad
		\varphi\bigl(U(x+Rx)\bigr)=-2a_0,
		\]
		so $\psi(x+Rx)=0$. It follows that for every $z\in M_0$,
		\[
		T^2(z,0)=(U^2z,0),
		\qquad
		(T^2)^*(z,0)=(U^{*2}z,0).
		\]
		Since $M_0$ reduces $U^2$, the nonzero subspace
		\[
		M=M_0\oplus\{0\}\subset H_0\oplus\mathbb C
		\]
		reduces $T^2$.
		
		The map
		\[
		J:H_{\mathrm e}\to M_0,
		\qquad
		Jx=\frac{x+Rx}{\sqrt5},
		\]
		is unitary, because $H_{\mathrm e}\perp H_{\mathrm o}$ and
		$\|Rx\|=2\|x\|$. Moreover, $J V_{\mathrm e}=U^2J$. Hence
		$U^2|_{M_0}$ is unitarily equivalent to the bilateral shift
		$V_{\mathrm e}$. In particular,
		\[
		\omega(T^2|_M)=\omega(V_{\mathrm e})=1,
		\]
		but this numerical radius is not attained: if a unit vector $x$ satisfied
		$|\langle V_{\mathrm e}x,x\rangle|=1$, equality in Cauchy--Schwarz would
		force $V_{\mathrm e}x=\lambda x$ for some $|\lambda|=1$, whereas a bilateral
		shift has no eigenvalues. Thus the restriction of $T^2$ to the nonzero
		reducing subspace $M$ fails to attain its numerical radius, and
		$T^2\notin\theta(H)$.
	\end{example}

	\section{Topological properties and density}\label{sec:topology}

	In this section we study the topological structure of $\theta(H)$ with respect to the operator norm and to other natural topologies. The main result, Theorem~\ref{thm:norm-density}, shows that $\theta(H)$ is norm-dense in $\mathcal{B}(H)$. If $H$ is finite-dimensional, then $\theta(H)=\mathcal B(H)$ and the density assertion is immediate; the arguments involving $W_e$ below are for infinite-dimensional $H$. The example of nonclosedness is likewise an infinite-dimensional phenomenon. The present study is related in spirit to the study of the essential numerical range and to that of operators with closed numerical ranges \cite{Bourin2003,JiLiang2018}, although the class $\theta(H)$ is governed by numerical-radius attainment on reducing subspaces rather than by closedness of the numerical range.
	
	We first show that $\theta(H)$ is not norm closed in $\mathcal{B}(H)$
	when $H$ is infinite-dimensional.
	
	\begin{example}\label{ex:nonclosed}
		Assume that $H$ is infinite-dimensional and choose an orthonormal basis
		$(e_n)_{n\in\mathbb N}$ of $H$. Define $T$ as in
		Example~\ref{ex:infinite-sum}, namely
		\[
		Te_n=\left(1-\frac1n\right)e_n.
		\]
		Then, by Theorem~\ref{thm:normal}, $T\notin\theta(H)$, since the moduli $1-\frac1n$ are strictly increasing. For each $N\in\mathbb{N}$ define
		\[
		T_Ne_n=
		\begin{cases}
			(1-\frac1n)e_n,&n<N,\\
			e_n,&n\geqslant N.
		\end{cases}
		\]
		Then $T_N$ is normal and its eigenvalue-modulus set is
		\[
		\left\{1-\frac1n:n<N\right\}\cup\{1\},
		\]
		which is reverse well-ordered. Hence
		\[
		T_N\in\theta(H).
		\]
		Moreover,
		\[
		\|T_N-T\|
		=
		\sup_{n\geqslant N}\frac1n
		=
		\frac1N
		\to0.
		\]
		Thus $\theta(H)$ is not norm closed.
	\end{example}
	
	\begin{corollary}[Norm closedness]\label{cor:norm-closedness}
		For a separable complex Hilbert space $H$, the class $\theta(H)$ is
		norm closed if and only if $\dim H<\infty$.
	\end{corollary}
	
	\begin{proof}
		If $H$ is finite-dimensional, then $\theta(H)=\mathcal B(H)$.
		If $H$ is infinite-dimensional, choose an orthonormal basis indexed by
		$\mathbb N$ and apply Example~\ref{ex:nonclosed}.
	\end{proof}
	
	In view of the preceding example, it is natural to ask whether $\theta(H)$ is nevertheless dense in $\mathcal{B}(H)$. In fact, the answer is affirmative. The proof of density relies on two fundamental results. First, we show that the set where the numerical radius strictly exceeds the essential numerical radius is open and norm-dense; this is Lemma~\ref{lem:gap-open-dense}. Second, we invoke the classical theorem of Halmos \cite{Halmos1968}, which states that the set of irreducible operators is norm-dense in $\mathcal{B}(H)$. Combining these two facts with Proposition~\ref{prop:strict-gap-attainment}, we obtain the norm density of $\theta(H)$.
	
	\begin{lemma}\label{lem:gap-open-dense}
		Assume that $H$ is infinite-dimensional. The set
		\[
		\mathcal{G}:=\{A\in\mathcal{B}(H):\omega(A)>\omega_e(A)\}
		\]
		is open and dense in $\mathcal{B}(H)$ with respect to the operator norm.
	\end{lemma}
	
	\begin{proof}
		For $A,B\in\mathcal{B}(H)$, the definition of the numerical radius gives
		\[
		|\omega(A)-\omega(B)|\leqslant\|A-B\|.
		\]
		We also have
		\[
		|\omega_e(A)-\omega_e(B)|\leqslant\|A-B\|.
		\]
		To justify the second estimate, put $d=\|A-B\|$. If $z\in W_e(A)$,
		there is a weakly null sequence of unit vectors $(x_n)$ such that
		$\langle Ax_n,x_n\rangle\to z$. The sequence $\langle Bx_n,x_n\rangle$
		is bounded, so after passing to a subsequence it converges to some
		$w\in W_e(B)$. Since
		$|\langle(A-B)x_n,x_n\rangle|\leq d$, we obtain $|z-w|\leq d$.
		Interchanging $A$ and $B$ gives the reverse estimate. Hence the Hausdorff
		distance between $W_e(A)$ and $W_e(B)$ is at most $d$, which directly yields
		$|\omega_e(A)-\omega_e(B)|\leq d$. Consequently,
		$A\mapsto\omega(A)-\omega_e(A)$ is continuous, so $\mathcal{G}$ is open.
		
		To prove density, fix $T\in\mathcal{B}(H)$ and $\varepsilon>0$. If
		$\omega(T)>\omega_e(T)$, then $T\in\mathcal{G}$. Otherwise, since
		$\omega_e(T)\leqslant\omega(T)$, we have
		$\omega(T)=\omega_e(T)=r$.
		
		If $r=0$, then $T=0$, because the numerical radius is a norm. Choose a
		nonzero rank-one operator $F$ with $\|F\|<\varepsilon$. Since $F$ is
		compact, its essential numerical range is $W_e(F)=\{0\}$, so
		$\omega_e(F)=0<\omega(F)$. Hence $F\in\mathcal{G}$.
		
		Suppose $r>0$. Choose $0<\delta<\min\{\varepsilon,2r\}$ and a unit
		vector $x$ such that
		\[
		|\langle Tx,x\rangle|>r-\frac{\delta}{2}.
		\]
		Set $a=\langle Tx,x\rangle$, so $a\ne0$, and let $P_x$ be the orthogonal
		projection onto $\mathbb{C}x$. Choose $\gamma\in\mathbb{C}$ with
		$|\gamma|=\delta$ so that
		$\langle\gamma P_xx,x\rangle=\frac{\delta a}{|a|}$, and put $F=\gamma P_x$.
		Then $F$ has rank one, $\|F\|=\delta$, and
		\[
		\omega(T+F)\geqslant|\langle(T+F)x,x\rangle|
		=|a|+\delta>r+\frac{\delta}{2}>r.
		\]
		Compact-perturbation invariance gives
		$\omega_e(T+F)=\omega_e(T)=r$. Therefore $T+F\in\mathcal{G}$ and
		$\|T-(T+F)\|=\delta<\varepsilon$. This proves density.
	\end{proof}
	
	\begin{remark}\label{rem:gap-vs-closed-range}
		Lemma~\ref{lem:gap-open-dense} is closely related to two earlier results.
		Bourin \cite[Proposition~1.1]{Bourin2003} recorded that every operator $A$
		admits a compact perturbation $A+K$ with $W_e(A)=\overline{W(A+K)}$, and Berg and Sims
		\cite{BergSims1984} proved that the numerical-radius attaining operators are
		norm dense. The gap condition $\omega(A)>\omega_e(A)$ used above is a
		quantitative form of the same circle of ideas: it guarantees attainment
		(Proposition~\ref{prop:strict-gap-attainment}), and the set where it holds is
		open and dense.
	\end{remark}
	
	We can now combine Lemma~\ref{lem:gap-open-dense} with Halmos's theorem to prove the norm density of $\theta(H)$.
	
	\begin{theorem}[Norm density]\label{thm:norm-density}
		If $H$ is a separable complex Hilbert space, then
		\[
		\overline{\theta(H)}^{\|\cdot\|}=\mathcal{B}(H).
		\]
	\end{theorem}
	
	\begin{proof}
		If $H$ is finite-dimensional, then $\theta(H)=\mathcal B(H)$ and the
		claim is immediate. Assume that $H$ is infinite-dimensional. Let
		$T\in\mathcal{B}(H)$ and $\varepsilon>0$. By
		Lemma~\ref{lem:gap-open-dense}, choose $A\in\mathcal{G}$ such that
		\[
		\|T-A\|<\frac{\varepsilon}{2}.
		\]
		Since $\mathcal{G}$ is open, there exists $\rho>0$ such that
		\[
		\{B\in\mathcal{B}(H):\|B-A\|<\rho\}\subseteq\mathcal{G}.
		\]
		By the classical density theorem of Halmos \cite{Halmos1968}, the set of irreducible operators is norm-dense in $\mathcal{B}(H)$. Thus we can choose an irreducible operator $R$ such that
		\[
		\|R-A\|<\min\left\{\frac{\rho}{2},\frac{\varepsilon}{2}\right\}.
		\]
		Then $R\in\mathcal{G}$, so
		\[
		\omega(R)>\omega_e(R).
		\]
		Proposition~\ref{prop:strict-gap-attainment} shows that $R$ attains its
		numerical radius. Since $R$ is irreducible, its only nonzero closed reducing
		subspace is $H$ itself. Hence $R\in\theta(H)$. Finally,
		\[
		\|T-R\|
		\leqslant\|T-A\|+\|A-R\|
		<\varepsilon.
		\]
		As $T$ and $\varepsilon$ were arbitrary, the theorem follows.
	\end{proof}
	
	\begin{remark}\label{rem:berg-sims}
		Theorem~\ref{thm:norm-density} strengthens the density theorem of Berg and Sims
		\cite{BergSims1984} in the property imposed on the approximating operators. Every operator
		in $\theta(H)$ attains its numerical radius, so the norm density of
		$\theta(H)$ implies the norm density of the numerical-radius attaining
		operators; moreover, the approximating operators can be chosen to satisfy
		the stronger reducing-subspace attainment property. The proof above,
		however, uses a different mechanism: instead of producing the approximant as a
		single compact perturbation of the original operator, it combines the open dense gap set $\{\omega>\omega_e\}$ with
		Halmos's density theorem for irreducible operators.
	\end{remark}
	
	As an immediate consequence of norm density and of the fact that norm convergence implies convergence in the strong, strong-star and weak operator topologies, we obtain the corresponding density statements in the weaker topologies.
	
	\begin{corollary}\label{cor:weak-density}
		If $H$ is a separable complex Hilbert space, then
		\[
		\overline{\theta(H)}^{\mathrm{SOT}}
		=\overline{\theta(H)}^{\mathrm{SOT}^*}
		=\overline{\theta(H)}^{\mathrm{WOT}}
		=\mathcal{B}(H).
		\]
	\end{corollary}
	
	\begin{proof}
		The conclusion follows immediately from Theorem~\ref{thm:norm-density},
		since convergence in operator norm implies strong, strong-star and weak
		operator convergence.
	\end{proof}
	
	\section{Concluding remarks}
	
	The results of this paper reveal a clear distinction between the normal and
	non-normal parts of the class $\theta(H)$. For normal operators, membership
	admits a complete spectral characterization in terms of pure point spectrum
	and the reverse well-ordering of eigenvalue moduli. Without normality, this
	spectral description no longer applies: compact and irreducible examples show
	that pure point spectrum need not be necessary, and the square of a non-normal
	operator in $\theta(H)$ may fail to belong to $\theta(H)$.
	
	The algebraic and perturbation results also exhibit different forms of
	stability. The class is preserved by finite orthogonal direct sums, but it is
	not generally preserved by translations, infinite orthogonal direct sums, or
	positive powers of non-normal operators. Moreover, although every compact
	operator belongs to $\theta(H)$, a finite-rank perturbation can destroy
	membership. Theorem~\ref{thm:invariant} gives the exact boundary of uniform
	stability: an operator remains in $\theta(H)$ under every compact perturbation
	if and only if it is compact.
	
	Finally, the essential numerical radius links the perturbation theory to the
	topological results. The strict gap $\omega(A)>\omega_e(A)$ guarantees
	numerical-radius attainment, and its open dense locus combines with Halmos's
	density theorem for irreducible operators to yield norm density of
	$\theta(H)$. Thus $\theta(H)$ is norm-dense, and consequently dense in the
	strong, strong-star, and weak operator topologies, while it is not norm closed
	in the infinite-dimensional setting. In finite dimensions,
	$\theta(H)=\mathcal B(H)$.
	
	\section{Further questions}
	
	Theorem~\ref{thm:general} gives an exact eigenvalue criterion for arbitrary
	operators, but it requires testing every nonzero reducing subspace. A natural
	problem is to find structural criteria for non-normal operators that avoid
	this quantification, for example in terms of the commutant $\{T,T^*\}'$.
	
	A second question concerns a fixed noncompact operator $T\in\theta(H)$:
	when does there exist $\varepsilon>0$ such that
	\[
	T+K\in\theta(H)
	\quad\text{for every }K\in\mathcal K(H)\text{ with }\|K\|<\varepsilon?
	\]
	Theorem~\ref{thm:invariant} addresses perturbations of arbitrary size.
	The scaled examples in Proposition~\ref{prop:rankone-instability} vary
	both $T$ and the perturbation, so they do not settle this local question
	for a fixed $T$.
	
	\section*{Acknowledgments}
	The authors acknowledge the use of DeepSeek V4.1 solely for literature search and language polishing. DeepSeek V4.1 was not used to generate mathematical content, proofs, or substantive research ideas. The authors are solely responsible for the content and accuracy of this paper. The authors also thank Professor Youqing Ji for valuable comments and suggestions.
	
	\section*{Declarations}
	
	\noindent \textbf{Author Contributions}
	All authors contributed equally.
	
	\bigskip
	
	\noindent\textbf{Funding}
	This work was supported by the National Natural Science Foundation of China (Grant No. 12471123).
	
	\bigskip
	
	\noindent \textbf{Data Availability}
	No datasets were generated or analyzed during the current study.

	\bigskip

	\noindent \textbf{Competing Interests}
	The authors declare no competing interests.
	
	\appendix
	
	\section{Closed operators and graph-subspace tools}\label{app:unbounded}
	
	This appendix collects the closed-operator facts and graph-subspace
	lemmas used in Section~\ref{sec:finite-direct-sum}. We first recall two
	standard facts; see \cite[Chapter X]{Conway1990} and \cite{kato1976}.
	We then give self-contained proofs of
	Lemmas~\ref{lem:graph-completion}--\ref{lem:graph-bounded-extension}.
	
	\subsection{Adjoints and positive square roots}
	Let $S:D(S)\subset H_1\to H_2$ be a closed linear operator, and put
	$K_1=\overline{D(S)}$. We regard $S$ as a densely defined operator from
	$K_1$ into $H_2$.
	
	\begin{enumerate}
		\item[(U1)] \emph{Adjoint of a closed operator.} The adjoint
		$S^*:D(S^*)\subset H_2\to K_1$ is closed and densely defined,
		and $S^{**}=S$. Moreover,
		\[
		\langle Sx,y\rangle=\langle x,S^*y\rangle
		\qquad(x\in D(S),\ y\in D(S^*)).
		\]
		
		\item[(U2)] \emph{Positive square root.} The operator $W=S^*S$ is
		positive and self-adjoint on $K_1$ with domain
		\[
		D(W)=\{x\in D(S):Sx\in D(S^*)\}.
		\]
		If $S$ is injective, then $W$ is injective. The positive square root
		$W^{\frac12}$ is a positive self-adjoint operator with domain $D(W^{\frac12})=D(S)$,
		and
		\[
		\|W^{\frac12}x\|=\|Sx\|\qquad(x\in D(S)).
		\]
	\end{enumerate}
	
	\subsection{Graph-norm and intertwining lemmas}\label{app:graph-tools}
	
	We begin with the graph-norm completion needed in the decomposition argument.
	
	\begin{lemma}[Graph-norm completion]\label{lem:graph-completion}
		Let $S:D(S)\subset H_1\to H_2$ be a closed operator whose domain is dense in a
		closed subspace $K_1\subset H_1$. Let $\widetilde K_1$ be the completion of
		$D(S)$ with respect to the graph norm $\|\cdot\|_S$. Then the graph map
		\[
		U_0:D(S)\to \operatorname{graph}(S),\qquad U_0x=(x,Sx),
		\]
		extends to a unitary operator from $\widetilde K_1$ onto
		$\operatorname{graph}(S)$.
	\end{lemma}
	
	\begin{proof}
		Since $S$ is closed, $\operatorname{graph}(S)$ is a closed subspace of
		$H_1\oplus H_2$. For $x,y\in D(S)$,
		\[
		\langle U_0x,U_0y\rangle_{H_1\oplus H_2}
		=\langle x,y\rangle+\langle Sx,Sy\rangle
		=\langle x,y\rangle_S.
		\]
		Thus $U_0$ is isometric. Its range is exactly $\operatorname{graph}(S)$. Hence
		$U_0$ extends uniquely to a unitary from the graph-norm completion
		$\widetilde K_1$ onto $\operatorname{graph}(S)$.
	\end{proof}
	
	The next lemma proves the polar decomposition and intertwining relations
	in the form used in Section~\ref{sec:finite-direct-sum}.
	
	\begin{lemma}[Polar decomposition and intertwining]\label{lem:polar-intertwining}
		Let $S:D(S)\subset K_1\to L$ be a closed injective operator with dense
		domain in the Hilbert space $K_1$ and dense range in the Hilbert space $L$.
		Put $W=S^*S$. Then $W$ is positive self-adjoint and injective, and
		\[
		S=UW^{\frac12}
		\]
		for a unitary operator $U:K_1\to L$. Suppose that
		$A\in\mathcal B(K_1)$ and $B\in\mathcal B(L)$ satisfy
		\[
		SAx=BSx,\qquad SA^*x=B^*Sx
		\qquad(x\in D(S)).
		\]
		Then
		\[
		WA=AW\quad\text{on }D(W),\qquad
		W^{\frac12}A=AW^{\frac12}\quad\text{on }D(S),
		\]
		and
		\[
		UA=BU.
		\]
		Consequently, $B=UAU^*$.
	\end{lemma}
	
	\begin{proof}
		Since $S$ is closed and densely defined on $K_1$, the operator
		$W=S^*S$ is positive and self-adjoint on $K_1$. It is injective because
		$S$ is injective. The polar decomposition for closed operators gives
		$S=UW^{\frac12}$, where $U:K_1\to L$ is unitary.
		
		The intertwining relations imply that $A(D(S))\subset D(S)$ and
		$A^*(D(S))\subset D(S)$. For $y\in D(S^*)$ and $x\in D(S)$,
		we have
		\[
		\begin{aligned}
			\langle Sx,By\rangle
			&=\langle B^*Sx,y\rangle
			=\langle SA^*x,y\rangle\\
			&=\langle A^*x,S^*y\rangle
			=\langle x,AS^*y\rangle.
		\end{aligned}
		\]
		Since $D(S)$ is dense in $K_1$ and $A$ is bounded, the definition
		of the adjoint gives $By\in D(S^*)$ and
		\[
		S^*By=AS^*y.
		\]
		Thus $B(D(S^*))\subset D(S^*)$ and $S^*B=AS^*$ on $D(S^*)$.
		
		We now check the domain of the commutation relation for $W$. Let
		$x\in D(W)$, so $x\in D(S)$ and $Sx\in D(S^*)$. By the
		domain invariance $B(D(S^*))\subset D(S^*)$ just proved,
		$BSx\in D(S^*)$. Since $Ax\in D(S)$ and $SAx=BSx$, we have
		$SAx\in D(S^*)$ and hence $Ax\in D(W)$. Moreover,
		\[
		\begin{aligned}
			WAx
			&=S^*SAx
			=S^*BSx\\
			&=AS^*Sx
			=AWx.
		\end{aligned}
		\]
		Therefore $AD(W)\subset D(W)$ and $WAx=AWx$ on $D(W)$.
		
		For $z\in\mathbb C\setminus\mathbb R$, the resolvent
		$(W-zI)^{-1}$ is a bounded operator on $K_1$. If $y\in K_1$ and
		$x=(W-zI)^{-1}y$, then $x\in D(W)$ and
		\[
		(W-zI)Ax=A(W-zI)x=Ay.
		\]
		Applying $(W-zI)^{-1}$ gives
		\[
		A(W-zI)^{-1}y=(W-zI)^{-1}Ay.
		\]
		Thus $A$ commutes with the resolvent of $W$, hence with the spectral
		projections of $W$ and its bounded Borel functional calculus. In
		particular, $A$ preserves $D(W^{\frac12})=D(S)$ and
		\[
		W^{\frac12}Ax=AW^{\frac12}x
		\qquad(x\in D(S)).
		\]
		
		Finally, for $x\in D(S)$,
		\[
		\begin{aligned}
			UAW^{\frac12}x
			&=UW^{\frac12}Ax
			=SAx\\
			&=BSx
			=BUW^{\frac12}x.
		\end{aligned}
		\]
		The range of $W^{\frac12}$ is dense in $K_1$, because $W^{\frac12}$ is
		self-adjoint and injective. Both $UA$ and $BU$ are bounded, so
		$UA=BU$ on all of $K_1$. Consequently $B=UAU^*$.
	\end{proof}
	
	Finally, we record the graph-norm boundedness of the self-adjoint parts
	$H_\zeta(A)$; this is the only additional unbounded-operator fact needed in
	the proof of the finite direct-sum theorem.
	
	\begin{lemma}[Bounded extension in the graph norm]\label{lem:graph-bounded-extension}
		Let $S:D(S)\subset K_1\to L$, $A\in\mathcal B(K_1)$, and
		$B\in\mathcal B(L)$ satisfy the assumptions of
		Lemma~\ref{lem:polar-intertwining}. For $\zeta\in\mathbb{T}$, put
		\[
		H_\zeta(A)=\frac{\zeta A+\overline{\zeta}A^*}{2},
		\qquad
		H_\zeta(B)=\frac{\zeta B+\overline{\zeta}B^*}{2}.
		\]
		Then $H_\zeta(A)$ is symmetric on $D(S)$ with respect to the graph inner
		product and bounded with respect to the
		graph norm $\|\cdot\|_S$. Consequently, it extends uniquely to a bounded
		self-adjoint operator $\widetilde H_\zeta$ on the graph-norm completion
		$\widetilde K_1$. Moreover,
		\[
		H_\zeta\bigl((A\oplus B)|_{\operatorname{graph}(S)}\bigr)
		=U_0\widetilde H_\zeta U_0^{-1}.
		\]
	\end{lemma}
	
	\begin{proof}
		The two intertwining relations imply that both $A$ and $A^*$ map
		$D(S)$ into itself. Consequently $H_\zeta(A)(D(S))\subset D(S)$,
		and, for $x\in D(S)$,
		\[
		SH_\zeta(A)x=H_\zeta(B)Sx.
		\]
		It follows that
		\[
		\begin{aligned}
			\|H_\zeta(A)x\|_S^2
			&=\|H_\zeta(A)x\|^2
			+\|H_\zeta(B)Sx\|^2\\
			&\leqslant C_\zeta\bigl(\|x\|^2+\|Sx\|^2\bigr)
			=C_\zeta\|x\|_S^2,
		\end{aligned}
		\]
		where
		$C_\zeta=\max\{\|H_\zeta(A)\|^2,\|H_\zeta(B)\|^2\}$.
		Thus $H_\zeta(A)$ has a unique bounded extension
		$\widetilde H_\zeta$ to the graph-norm completion $\widetilde K_1$.
		
		We next verify symmetry in the graph inner product. For $x,y\in D(S)$,
		\[
		\begin{aligned}
			\langle H_\zeta(A)x,y\rangle_S
			&=\langle H_\zeta(A)x,y\rangle
			+\langle SH_\zeta(A)x,Sy\rangle\\
			&=\langle H_\zeta(A)x,y\rangle
			+\langle H_\zeta(B)Sx,Sy\rangle\\
			&=\langle x,H_\zeta(A)y\rangle
			+\langle Sx,H_\zeta(B)Sy\rangle\\
			&=\langle x,H_\zeta(A)y\rangle_S.
		\end{aligned}
		\]
		Here we used the self-adjointness of $H_\zeta(A)$ and
		$H_\zeta(B)$ in their respective Hilbert spaces. By continuity, the
		symmetry identity extends to all of $\widetilde K_1$. A bounded symmetric
		operator defined on the whole Hilbert space is self-adjoint, so
		$\widetilde H_\zeta$ is self-adjoint.
		
		Finally, the graph map $U_0x=(x,Sx)$ intertwines the two self-adjoint
		parts on $D(S)$; by continuity this gives
		\[
		H_\zeta\bigl((A\oplus B)|_{\operatorname{graph}(S)}\bigr)
		=U_0\widetilde H_\zeta U_0^{-1}.
		\]
	\end{proof}

\end{document}